\documentclass[11pt,a4paper]{article}

\usepackage{titlesec}

\usepackage[utf8]{inputenc}

\usepackage[T1]{fontenc}
\usepackage[english]{babel}
\usepackage{amsmath}
\usepackage{amssymb}
\usepackage{amsthm}
\usepackage{dsfont}
\usepackage{constants}

\usepackage{graphicx, psfrag}
\usepackage{xcolor}
\usepackage{subfigure}

\usepackage{mathtools}
\usepackage[titletoc,toc,title]{appendix}

\usepackage{color}
\usepackage{caption}
\usepackage{bbm}
\usepackage[square, numbers, sort]{natbib}
\usepackage{float}
\usepackage{multirow}
\usepackage{array}
\usepackage{geometry}
\usepackage{fancyhdr}
\usepackage[hyphens]{url}
\usepackage[pdfpagelabels=false]{hyperref}
\usepackage{lettrine}
\usepackage{theoremref}
\usepackage[]{todonotes}

\usepackage[normalem]{ulem}

\usepackage[shortlabels]{enumitem}

\numberwithin{equation}{section}

\setenumerate[0]{label=(\alph*)}

\newcommand{\1}[1]{\mathds{1}_{ \{ #1 \} }}
\newcommand{\floor}[1]{{\lfloor #1 \rfloor}}
\newcommand{\Z}{\ensuremath{\mathbb{Z}}}
\newcommand{\R}{\ensuremath{\mathbb{R}}}
\newcommand{\N}{\ensuremath{\mathbb{N}}}

\newcommand{\abs}[1]{\lvert #1 \rvert}

\renewcommand{\P}{\ensuremath{\mathbb{P}}}

\newcommand{\Pb}[2]{\mathbb{P}_{#1}\bigl( #2 \bigr)}
\newcommand{\Pbe}[3]{\mathbb{P}_{#1}^{#2}\bigl( #3 \bigr)}

\newcommand{\E}[2]{\mathbb{E}_{#1}\!\Bigl[ #2 \Bigr]}
\newcommand{\Ee}[3]{\mathbb{E}_{#1}^{#2}\Bigl[ #3 \Bigr]}
\newcommand{\Eeb}[3]{\mathbb{E}_{#1}^{#2}\bigl[ #3 \bigr]}

\newconstantfamily{c}{symbol=C}
\newconstantfamily{c_s}{symbol=c}

\newtheorem{theorem}{Theorem}[section]

\newtheorem{proposition}[theorem]{Proposition}
\newtheorem{lemma}[theorem]{Lemma}

\newtheorem{corollary}[theorem]{Corollary}

\newtheorem{remark}[theorem]{Remark}

\begin{document}

\title{Sublinear Power-Law Growth of Level Set Diameters in Heterogeneous F-KPP Equations}

\author{
  Alexander Drewitz
  \thanks{Universität zu Köln,
    Department Mathematik/Informatik,
    Weyertal 86--90,
    50931 Köln, Germany.
    Email: \url{adrewitz@uni-koeln.de}}
    \and
  	Pascal Hanigk
  	\thanks{Institut für Mathematik, 
    Johannes Gutenberg-Universität Mainz,
    Staudingerweg 9,
    55099 Mainz, Germany.
    Email: \url{pahanigk@uni-mainz.de}}
    \and
     Lisa Hartung
    \thanks{Institut für Mathematik, 
    Johannes Gutenberg-Universität Mainz,
    Staudingerweg 9,
    55099 Mainz, Germany \& RMU Co-Affiliate of Technical University Darmstadt, Germany
    Email: \url{lhartung@uni-mainz.de}}
}

\date{\today}

\maketitle

\begingroup
\renewcommand{\thefootnote}{}
\endgroup

\begin{abstract}
For every $\alpha\in(1/4,1)$, we construct a heterogeneous medium $\xi$ such that the diameter of the level sets of the solution to the F-KPP equation is of order $t^\alpha$ along a sequence of times tending to infinity.
This is a natural continuation of a series of recent works establishing logarithmic lower bounds as well as linear asymptotics for the diameters of the respective level sets in random and deterministic heterogeneous media.
Both results assume $\sup\xi > 2\inf\xi$; we also show that the diameter can diverge in the critical regime $\sup\xi=2\inf\xi$.

A key ingredient is an estimate for branching Brownian motion with fixed time horizon $t$ and branching rate $1$ on the negative half-line and $2+t^{-\eta}$ on the positive half-line.
This estimate exhibits a saturation transition: for \(\eta<1/2\) the polynomial and logarithmic corrections to the extremal position depend on $\eta$, whereas for $\eta\geq1/2$ the extremal position becomes independent of $\eta$, up to an additive constant. 
\end{abstract}

\tableofcontents

\thispagestyle{empty}

\section{Introduction}
\subsection{Discussion and previous results}

The classical F-KPP equation
\begin{align*}
    u_t(t,x) &= \frac12 u_{xx}(t,x) +u(t,x)(1-u(t,x)), \qquad (t,x)\in(0,\infty)\times\R, \\
    u(0,\cdot) &= u_0,
\end{align*}
with initial condition $u_0 : \R \to [0,1]$,
was introduced in the seminal works of Fisher \cite{fisher1937} and of Kolmogorov, Petrovskii and Piskunov \cite{kolmogorov1937}. 
It is well established (see \cite{kolmogorov1937, uchiyama77, bramson1983convergence}) that the classical F-KPP equation admits traveling wave solutions of every speed $\lambda\geq\sqrt{2}$, and for suitable initial conditions $u_0$ (e.g.\ Heaviside) the solution converges, up to a logarithmic shift, to the wave of minimal speed $\sqrt{2}$. 
In particular, the transition region, that is the region between the stable state $u=1$ and the unstable state $u=0$, remains bounded.
A natural generalization of the above equation is 
\begin{equation}\label{eq:KPP}\tag{F-KPP}
\begin{split}
    u_t(t,x) &= \frac{1}{2}u_{xx}(t,x) +  \xi(x)\, F(u(t,x)),
    \qquad (t,x)\in(0,\infty)\times\R,\\
    u(0,\cdot) &= u_0,
\end{split}
\end{equation}
with a more general nonlinearity $F$ and $\xi : \R\to(0,\infty)$, a function modeling a heterogeneous medium. 
In our specific setting, we consider nonlinearities $F$ that are generated by the probability generating function belonging to branching Brownian motion, see \eqref{eq:standard_cond} below.

In this article we are interested in the effects that a heterogeneous medium $\xi$ has on the geometry of the transition region of the solution.
For this purpose, we introduce, for the solution $u$ to \eqref{eq:KPP}, $\varepsilon\in(0,1)$ and $t>0$, the quantities 
\begin{equation}\label{def:front_KPP}
  \begin{split}
    m^{\varepsilon, +}(t) &
    := \sup\{ x \in \R \, : \, u(t,x) \geq \varepsilon\} \big(=\sup\{ x \in \R \, : \, u(t,x) = \varepsilon\}\big), \\
    m^{\varepsilon,-}(t) &
    :=  \inf\{ x \in \R \, : \,  u(t,x)\leq \varepsilon\} \big( =  \inf\{ x \in \R \, : \,  u(t,x)= \varepsilon\} \big), 
  \end{split}
\end{equation}
which represent the right and left $\varepsilon$-fronts, respectively.
The \emph{diameter of the $\varepsilon$-level set}
\begin{equation}\label{def:front_diam}
    I^\varepsilon(t) \coloneqq m^{\varepsilon, +}(t) - m^{\varepsilon, -}(t) \geq 0
\end{equation}
quantifies the spatial width between the right and left $\varepsilon$-fronts.

For the classical F-KPP equation with a non-increasing (and non-constant) initial condition $u_0:\R \to [0,1]$, the solution $u(t,\cdot)$ is continuous and strictly decreasing on $\R$ for every $t>0$. Consequently $m^{\varepsilon, +}(t) = m^{\varepsilon, -}(t)$, implying $I^\varepsilon(t)=0$ for all $t>0$.

In the heterogeneous setting, the monotonicity of the solution to \eqref{eq:KPP} may fail due to the medium $\xi$, potentially leading to unbounded growth of $I^\varepsilon(t)$.
In this case, so called \emph{generalized transition fronts} for \eqref{eq:KPP}, i.e.\ entire solutions connecting $1$ to $0$ whose transition region has uniformly bounded width (cf. \cite[Definition 1.2]{hamel2021diameters}), do not exist. This follows from the proof of Corollary~1.3 in \cite{hamel2021diameters}, which shows that the existence of generalized transition fronts implies $\sup_{t>0} I^\varepsilon(t)<\infty$ for all $\varepsilon\in(0,1)$, at least for Heaviside initial condition.

A first non-existence result for generalized transition fronts in heterogeneous media was given by Nolen, Roquejoffre, Ryzhik and Zlatoš \cite{nolen_zlatos_12} for bounded media satisfying $\xi(x)\ge1$ and $\lim_{\abs{x}\to\infty}\xi(x)=1$. 
Letting $\mu_+ \coloneqq \sup_{x\in\R}\xi(x)$ and $\mu_-\coloneqq \inf_{x\in\R}\xi(x)$, their condition for non-existence requires that $\mu_+>2\mu_-$.
Conversely, Zlatoš \cite{zlatos12} showed that, if $\mu_+<2\mu_-$, then there exists a generalized transition front for \eqref{eq:KPP}, which establishes $\mu_+=2\mu_-$ as a \emph{critical} regime.
In \cite{hamel2021diameters}, Hamel and Nadin constructed a function $\xi$ consisting of increasingly long stretches alternating between the values $\mu_-$ and $\mu_+>2\mu_-$ and proved that the diameters of the level sets $I^\varepsilon(t_n)$ can be bounded from below and above by linear functions along a suitable sequence of times $(t_n)$, thus extending the non-existence result in \cite{nolen_zlatos_12} to a new class of heterogeneous nonlinearities. 
A related result had previously been obtained by Černý, Drewitz and Schmitz \cite{cerny2021unbounded}, who considered the case of a random function $\xi$, exhibiting sufficient stochastic independence. They showed that along a sequence of time-space points $(t_n, x_n)$, one has that $I^\varepsilon(t_n)$ diverges at least logarithmically in $x_n$. Their result also requires $\mu_+>2\mu_-$.

It is the main goal of this article to extend the above findings by understanding the sublinear regime for $I^\varepsilon(t)$. More precisely, for any $\alpha\in(1/4,1)$, we construct a medium $\xi$ such that $I^\varepsilon(t_n)$ can be of order $t_n^\alpha$ along a sequence of times $(t_n)$, cf.\ Theorem~\ref{thm:main} (the restriction to $\alpha>1/4$ stems from our construction of $\xi$,
see Remark \ref{rem:sat} also).
To achieve this, we first analyze the \emph{near-critical} regime, where $\mu_+=\mu_+(t)\downarrow 2\mu_-$ as $t\to\infty$, and prove the existence of a family of media $(\xi_t)_{t>0}$, indexed by the time horizon $t$, in which $I^\varepsilon(t)$ diverges with sublinear speed (cf. Proposition~\ref{prop:main_simple}). 
As part of our investigation, we also construct an example in the critical regime $\mu_+=2\mu_-$ in which $I^\varepsilon(t)$ diverges (cf.\ Theorem~\ref{thm:main_crit}), confirming $\mu_+<2\mu_-$ as the sharp threshold for the universal existence of generalized transition fronts.

Our proofs are entirely probabilistic and rely heavily on the duality between the heterogeneous F-KPP equation \eqref{eq:KPP} and branching Brownian motion with spatially dependent branching rates, see Proposition~\ref{prop:mckean1} below.
A large portion of this article is dedicated to understanding the extremal behavior of branching Brownian motion at time $t$, where particles at position $x$ branch at rate $1$ for $x<0$ and rate $\mu_+(t)\geq2\mu_-$ for $x\geq0$, in the (near-)critical regime where $\mu_+(t)\downarrow2\mu_-$ (cf. Theorem~\ref{thm:main_bbm}).
We then apply our results for branching Brownian motion to prove our results pertaining the level sets of the solution to \eqref{eq:KPP}. 

In order to describe the admissible nonlinearities for \eqref{eq:KPP} that we will cover, with the convention $\N = \{ 1, 2, 3, \ldots\}$, let $(p_k)_{k\in\N}$ be an arbitrary sequence of reals in $[0,1]$ and define $F:[0,1]\to[0,1]$ via
\begin{equation}\label{eq:standard_cond}\tag{PROB}
    F(u) \coloneqq (1-u) - \sum_{k\in\N} p_k (1-u)^k,
\end{equation}
where
\begin{equation}\label{eq:mom_offspr}
    \sum_{k\in\N} p_k = 1, \quad \sum_{k\in\N}kp_k=2, \quad \text{and} \quad m_2 \coloneqq \sum_{k\in\N} k^2p_k<\infty.
\end{equation}
In passing, we note that $F'(0)=1$. 
Choosing $p_2=1$ recovers the classical F-KPP equation.
The reason for considering this type of nonlinearity is its suitability for being investigated using techniques from branching processes.

\subsection{Main results}
\emph{Notation.} We write $f = \Theta(g)$ if $cg(x) \leq f(x) \leq Cg(x)$ for some constants $c,C > 0$ and all $x$ large enough from the intersection of the domains of definition of $f$ and $g$. Similarly, we write $f=O(g)$ if $\abs{f(x)} \leq C\abs{g(x)}$ for some constant $C>0$ and all $x$ large enough, and $f=o(g)$ if $\lim_{x\to\infty} f(x)/g(x) =0$.

Our first main result states that there exists a medium $\xi$ in which the diameter of the $\varepsilon$-level set of the solution to \eqref{eq:KPP} diverges with sublinear power-law speed.

\begin{theorem} \label{thm:main}
Let $F$ satisfy \eqref{eq:standard_cond} and $u_0 = \mathds{1}_{(-\infty, 0]}$.  Then, for all $\alpha \in (1/4,1)$ and all $\varepsilon$ sufficiently small, there exists a function $\xi = \xi_{\alpha,\varepsilon}: \R \to (0,\infty)$ and an increasing sequence of times $(t_n)_{n\in\N}$ with $\lim_{n \to \infty} t_n = \infty$ such that
\begin{align*}
	I^\varepsilon(t_n) = \Theta\!\left( t_n^\alpha \right),
\end{align*}
where $I^\varepsilon$ is defined in \eqref{def:front_diam}.
\end{theorem}

The condition $\alpha\in(1/4,1)$ stems from our construction of $\xi$. We first introduce a simplified time-dependent medium $\widetilde{\xi}$, and then construct the time-independent medium $\xi$ that reproduces the desired properties of $\varepsilon$-fronts in the medium $\widetilde{\xi}$ along a suitable sequence of times. 
We postpone the discussion of the behavior of the transition region of $\widetilde{\xi}$, as well as the natural emergence of the restriction $\alpha\in(1/4,1)$ in this setting, to Section \ref{sec:simplePot}. See in particular Proposition~\ref{prop:main_simple} and the discussion below.

Our next result shows that in the critical regime $\sup_{x\in\R}\xi(x) = 2 \inf_{x\in\R}\xi(x)$ it is possible that $I^\varepsilon(t)$ diverges and therefore no generalized transition front exists.
\begin{theorem}\label{thm:main_crit}
    Let $F$ satisfy \eqref{eq:standard_cond} and $u_0 = \mathds{1}_{(-\infty, 0]}$. For all $\varepsilon$ sufficiently small there exists a function $\xi=\xi_\varepsilon : \R\to(0,\infty)$ with $\sup_{x\in\R}\xi(x) = 2 \inf_{x\in\R}\xi(x)$, a constant $c>0$ and an increasing sequence of times $(t_n)_{n\in\N}$ with $\lim_{n\to\infty}t_n=\infty$ such that
    \begin{equation*}
        I^\varepsilon(t_n) \geq c\,t_n^{1/4} (\log t_n)^{1/2}
    \end{equation*}
    for all $n$ large enough.
\end{theorem}

\subsection{Connections to branching Brownian motion and asymptotics around the front}  \label{sec:conn}

Our proofs of Theorems~\ref{thm:main} and \ref{thm:main_crit} heavily take advantage of a probabilistic representation of the solution to \eqref{eq:KPP}.

We define a \emph{branching Brownian motion in medium $\xi$} (BBM$\xi$) as follows: For a locally integrable function $\xi : \R\to(0,\infty)$,  an offspring law $(p_k)_{k\in\N}$, and fixed $x\in\mathbb{R}$, consider an initial particle starting at a point $x$ and moving as a standard Brownian motion $(B_t)_{t\geq0}$ on $\mathbb{R}$. While at site $y\in \R$, the particle dies at rate $\xi(y)$. More precisely, for an exponentially distributed random variable $S$ with parameter one, independent of everything else, the first particle dies at time $\inf\big\{t\geq0: S<\int_0^t\xi(B_s)\,{\rm d} s \big\}$. 
When the initial particle dies, it gives birth to $k$ new particles with probability $p_k$, $k\in\N.$ Note that our choice of $(p_k)_{k\in\N}$ implies that no particles dies without offspring.
The new particle(s) start their evolution at the site where their parent particle had died, and they evolve independently of everything else and according to the same stochastic behavior as their parent. 
We denote the corresponding probability measure by $\mathbb{P}_x^\xi$ and write $\mathbb{E}_{x}^\xi$ for the respective expectation. If the medium is constant equal to $c>0$, we write $\mathbb{P}_x^c$ and $\mathbb{E}_{x}^c$ instead.
Any probability measure and expectation without a superscript corresponds to the Brownian motion.
By $N(t)$ we denote the set of particles alive at time $t$ and for $\nu\in N(t)$ we let 
\begin{equation*}
(X_s^\nu)_{s\in[0,t]}
\end{equation*}
be the spatial trajectory of the genealogy of ancestral particles of $\nu$ up to time $t$. Furthermore, for $t\in [0,\infty),$ by $M_t$ we denote the position of the leftmost particle of BBM$\xi$ in $N(t)$, i.e.  $M_t \coloneqq \min \left\{ X^\nu_t \,:\, \nu\in N(t) \right\}$.

We are now in a position  to formulate a representation of $u(t,x)$ which is sometimes referred to as the McKean representation. 
The respective result in the homogeneous setting was proven in \cite{branch_markov_I} and \cite{mckean}. 

\begin{proposition}[Proposition 2.1 in \cite{DreSchmi19}]\label{prop:mckean1}
	Let $u_0$ be a sufficiently nice initial condition, let $\xi$ be locally H\"{o}lder continuous and $F$ fulfill \eqref{eq:standard_cond} with $(p_k)_{k\in \N}$ satisfying \eqref{eq:mom_offspr}. Then
	the solution to \eqref{eq:KPP} is given by
	\begin{equation}\label{eq:McKean}
	\tag{McKean}
	u(t,x) = 1- \Ee{x}{\xi}{ \prod_{\nu\in N(t)} \big(1-u_0(X_t^\nu)\big) },\quad (t,x)\in [0,\infty) \times\R.
	\end{equation}
\end{proposition}
\begin{remark}
    In particular, \eqref{eq:McKean} holds for the Heaviside-type initial condition $u_0 = \mathds{1}_{(-\infty, 0]}$. This choice yields
	\begin{equation}\label{eq:mckMin}
		u(t,x) = 1 - \P_{x}^\xi\big( X_t^\nu > 0 \,\forall\, \nu\in N(t) \big) = \Pbe{x}{\xi}{ M_t \leq 0 }
	\end{equation}	 
	for all $(t,x)\in [0,\infty)\times\R$.
\end{remark}
In the remainder of the paper, when proving properties of the solution $u(t,x)$ of \eqref{eq:KPP}, we will mostly do so by considering $\Pbe{x}{\xi}{ M_t \leq 0 }$ in combination with \eqref{eq:mckMin}.

As part of constructing our medium $\xi$ for Theorem~\ref{thm:main}, we study branching Brownian motion in media of the form
\begin{equation}\label{def:xiHat}
    \widehat{\xi}(z) \coloneqq \widehat{\xi}_t(z) \coloneqq \begin{cases}
        1, & z\in(-\infty,0), \\
        \mu_+, & z\in[0,\infty),
    \end{cases}
\end{equation}
where 
\begin{equation} \label{def:mu+}
    \mu_+ = \mu_+(t) = 2+t^{-\eta}
\end{equation}
with $\eta=\eta(t)>0$ satisfying
\begin{equation}\label{def:eta}
    \eta=\eta(t) = \eta_0 + \varphi(t)
\end{equation}
for some fixed $\eta_0\in(0,\infty]$ (with the convention that $t^{-\infty}=0$ so that $\mu_+=2$ is allowed) and a function $\varphi(\cdot)\geq0$, such that $\lim_{t\to\infty}\varphi(t)=0$ and the function $t^{-\varphi(t)}$ is non-increasing in $t$. The choices of $\varphi$ used throughout this paper are
\[
    \varphi\equiv0 \quad\text{and}\quad \varphi(t) = \frac{\log(\log(t)\vee1)}{\log(t)}.
\]
We establish tightness of the re-centered leftmost particle at time $t$ for BBM$\widehat{\xi}$ started at the origin. 

\begin{theorem}\label{thm:main_bbm}
    For $\mu_+$ as in \eqref{def:mu+} define the function
    \begin{equation}\label{def:m(t)}
        \widehat{m}(t) \coloneqq \frac{\mu_+ t}{\sqrt{2(\mu_+ -1)}} - \frac{3}{2\sqrt{2}}\log(t) - \frac{1}{2\sqrt{2}}\log(t^{-\eta}\vee t^{-1/2}).
    \end{equation}
   Then, for any function $\eta=\eta(t)$ as in \eqref{def:eta}, the family of distributions
    \[
        \mathbb P_0^{\widehat{\xi}_t} \circ( M_t + \widehat{m}(t) )^{-1}, \qquad t\in[1,\infty)
    \]
    is tight.
\end{theorem}
\begin{remark}
\begin{enumerate}
    
\item
    A key observation is that the asymptotics of $\widehat{m}(t)$ become insensitive to the precise choice of $\mu_+$ in \eqref{def:mu+} provided that $\mu_+$ converges sufficiently fast.
    Indeed, Taylor's formula implies
	\begin{align}\label{eq:taylor}
		\frac{\mu_+ t}{\sqrt{2(\mu_+ -1)}} = \sqrt{2}t + \Theta\bigl( (\mu_+ -2)^2 t \bigr) = \sqrt{2}t + \Theta\bigl( t^{1-2\eta} \bigr).
	\end{align}
	For $\eta_0\in[1/2,\infty]$,
	\begin{align*}
		\widehat{m}(t) = \sqrt{2}t - \frac{3-\frac12}{2\sqrt{2}}\log(t) + O(1).
	\end{align*}
    For standard BBM (where $\xi\equiv1$), it is well known that $M_t + (\sqrt{2}t-\frac{3}{2\sqrt{2}}\log(t))$ is tight. 
	In particular, the position of $M_t$ for BBM$\widehat{\xi}$ remains smaller than that of standard BBM by $\log(t)/(4\sqrt{2})$. 

\item
The transition at $\eta_0=1/2$ arises because a particle attempting to become minimal first spends some time to the right of the origin, where the typical excursion time above $0$ turns out to be of order $t^{1-\eta}$ (see \eqref{def:r*}) with fluctuations of order $\sqrt{t}$. Consequently, when $\eta_0\in(1/2,\infty]$ (and similarly for $\eta_0=1/2$), these fluctuations dominate the typical excursion time, and this strategy reaches a form of equilibrium. More precisely, over time scales of order $\sqrt{t}$, the advantage of the favorable branching rate is balanced by the disadvantage of delaying the descent into the negative half-line. Since this time scale does not depend on the particular choice of $\mu_+$, the correct centering given by $\widehat{m}(t)$ defined in \eqref{def:m(t)} likewise becomes insensitive to the exact choice of $\eta$.
\end{enumerate}
\end{remark}
  
\subsection{Overview of the remainder of the article}
The rest of this paper is devoted to the proofs of Theorems~\ref{thm:main}, \ref{thm:main_crit} and \ref{thm:main_bbm}. In order to do so, we will ultimately construct an increasing sequence of times $(t_n)$ with $\lim_{n\to\infty} t_n = \infty$ and a corresponding medium such that the properties of Theorems~\ref{thm:main} and \ref{thm:main_crit} are fulfilled (see Theorem~\ref{thm:main_gen}).
Heuristically, $\xi$ (later defined in \eqref{def:xi}) will consist of increasingly long stretches where the medium alternates between $1$ and values $\mu_+$ slightly larger than or equal to $2$, see \eqref{def:mu+} for further details. A crucial role in proving Theorems \ref{thm:main} and \ref{thm:main_crit} will be played by understanding the solution around those locations in space where a transition from $1$ to $\mu_+$ of the medium occurs. Using Proposition~\ref{prop:mckean1}, this can be transferred to understanding the distribution of the leftmost particle of BBM$\xi$ started around such points at suitable times.

In order to facilitate the reader's understanding, we therefore start with understanding the leftmost particle of BBM$\widetilde{\xi}$ and BBM$\widehat{\xi}$, where $\widetilde{\xi}$ is a simplified medium when compared to $\xi$, while $\widehat{\xi}$ is a simplified and shifted medium when compared to $\xi$. These media will take the value $1$ to the left of a certain point $\widetilde{x}$ and the origin, respectively, and the value $\mu_+$ to the right of these points. See \eqref{def:xiTilde} and \eqref{def:xiHat} for the definitions of $\widetilde{\xi}$ and $\widehat{\xi}$, respectively.

Section \ref{sec:minSimpl} is dedicated to analyzing the leftmost particle position of BBM$\widehat{\xi}$ and the proof of Theorem~\ref{thm:main_bbm}. Furthermore we study the effects of varying the initial position near the origin, the interface of the medium $\widehat{\xi}$, which will be needed later in the proofs of Theorems~\ref{thm:main} and \ref{thm:main_crit}. 
The results are stated in Section~\ref{sec:Result}.

Section \ref{sec:simplePot} is dedicated to understanding the $\varepsilon$-fronts in the simplified medium $\widetilde{\xi}$. In order to do so we start by defining and analyzing the center point $\widetilde{x}$ and then consider starting points close to $\widetilde{x}$. This can be deduced from the results of Section~\ref{sec:minSimpl}.

In Section \ref{sec:fullPot} we then define the medium $\xi$ and show how to transfer the results from the simplified medium $\widetilde{\xi}$ to the full setting $\xi$. Theorems~\ref{thm:main} and \ref{thm:main_crit} are immediate consequences of Theorem~\ref{thm:main_gen}.

We use numbered constants, such as $C_1, C_2, \ldots$, to represent positive finite constants appearing at various points throughout the article. Within a proof, generic constants, denoted by $c$, $C$, and similar symbols, may change their value from one line to the next within an inequality.

\bigskip

\noindent\textbf{Statement on AI usage.} 
We used Opus 5.5 for proofreading and to check the text for grammatical, notational and stylistic inconsistencies.
In addition, the elementary bound \eqref{eq:sup_gauss} was suggested by AI, which simplified the proofs of Lemma~\ref{lem:E_X_int2.0} and Proposition~\ref{prop:E_X}.
AI tools were not used to develop mathematical ideas or in the writing of the first draft.
All suggested changes were reviewed by the authors, who take full responsibility for the content
of the paper.

\bigskip

\noindent\textbf{Acknowledgment.} The work of PH and LH was partly funded by the Deutsche Forschungsgemeinschaft (DFG, German Research Foundation) through Project-ID 44391315 within SPP 2265.

\section{Bounding the minimal particle position of BBM$\widehat{\xi}$} \label{sec:minSimpl}
\begin{figure}[htbp]
\centering
\includegraphics[width=0.8\textwidth]{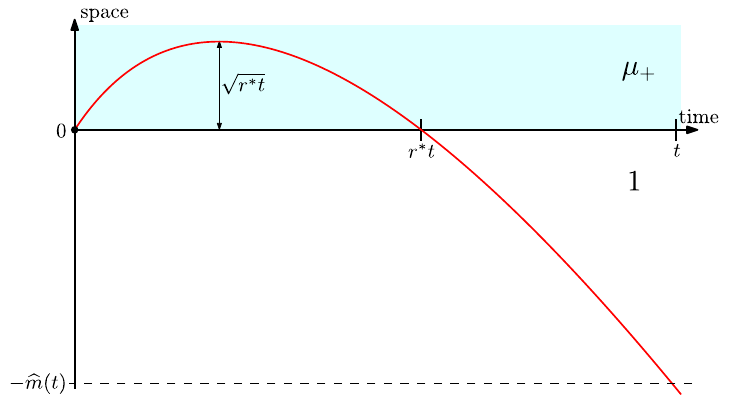}
\caption{A sketch of the typical path of the leftmost particle of BBM$\widehat{\xi}$ at time $t$ starting at the origin; for spatial values in the turquoise area (corresponding to real numbers in the positive half-line) the medium takes value $\mu_+$, while in the white area (corresponding to real numbers in the negative half-line) it takes the value $1$.}
\label{fig:optimal_path}
\end{figure}

In this section we provide a proof of Theorem~\ref{thm:main_bbm}, and for that purpose we show that
\begin{equation*}
    \Pbe{0}{\widehat{\xi}}{M_t + \widehat{m}(t) \leq y}
\end{equation*}
is bounded away from $0$ and $1$ asymptotically, for appropriate choices of $y$. With regard to later estimating the $\varepsilon$-fronts, we also establish estimates for varying starting positions of BBM$\widehat \xi$. 
Define
\begin{equation}\label{def:r*} r^* \coloneqq \frac{\mu_+-2}{2(\mu_+-1)} = \frac{t^{-\eta}}{2(1+t^{-\eta})}\sim t^{-\eta}/2;
\end{equation}
we will see below (cf.\ \eqref{eq:IDef} and \eqref{eq:Z'}) that  (at least for $\eta_0\in(0,1/2)$) the value $r^*$ can be interpreted as the typical fraction of time a particle will spend above $0$ in order to reach an extreme particle position at time $t$ under $\P_{0}^{\widehat{\xi}}$
up to the accuracy we are interested in here (see Figure~\ref{fig:optimal_path}).

\subsection{Results} \label{sec:Result}
We distinguish  two cases, depending on how fast $\mu_+$ converges to $2$. First, we consider the case of slow convergence, which in our case means $\eta_0\in(0,1/2)$. Note that in this case \eqref{def:m(t)} evaluates to
\begin{align*}
	\widehat{m}(t) =\frac{\mu_+ t}{\sqrt{2(\mu_+ -1)}} - \frac{3-\eta}{2\sqrt{2}}\log(t)
\end{align*}
for $t$ large enough.

\begin{theorem}\label{the:min_case1}
    There exist constants $\Cl[c]{cst:ub1},\Cl[c]{cst:lb11} > 0$ such that for every $\eta_0\in(0,1/2)$, if $z=z(t)\geq0$ and $y=y(t)=o(\log t)$, then
    \begin{equation} \label{eq:etaSmallUB}
        \Pbe{z}{\widehat{\xi}}{M_t +\widehat{m}(t) \leq y} \leq \Cr{cst:ub1} (z\vee1)\,{\rm e}^{\sqrt{2}y-\frac{z^2}{2r^*t\,(1+o(1))}}\, +o(1)
    \end{equation}
    for all $t$ sufficiently large, where $o(1)$ denotes a function converging to $0$ as $t\to\infty$, independent of the choice of $z$.
    
    If also $z\leq (r^*t)^{3/4}$, then 
\begin{equation}\label{eq:etaSmallLB}
		\Pbe{z}{\widehat{\xi}}{M_t + \widehat{m}(t) \leq y} \geq \Cr{cst:lb11} \bigl( 1 - \exp\bigl\{-\Cr{cst:lb11} (z\vee1)\, {\rm e}^{\sqrt{2}y-\frac{z^2}{2r^*t\,(1+o(1))}}\bigr\} \bigr)
	\end{equation}
	for all $t$ large enough.
\end{theorem}

Our second result in this section deals with the case where the convergence rate of $\mu_+$ is fast in the sense that $\eta_0\in[1/2,\infty]$. Note that in this case, 
\begin{align*}
	\widehat{m}(t) = \frac{\mu_+ t}{\sqrt{2(\mu_+ -1)}} - \frac{3-\frac12}{2\sqrt{2}}\log(t)
\end{align*}
for $t$ large enough.
    
\begin{theorem}\label{the:min_case2}
     Let $\eta_0\in[1/2,\infty]$.
    \begin{enumerate}[label=(\alph*)]
        \item There exists a constant $\Cl[c]{cst:ub2}>0$ such that if $z=z(t)\geq0$ and $y=y(t)=o(\log t)$, then
        \begin{equation} \label{eq:etaLargeUB}
            \Pbe{z}{\widehat{\xi}}{M_t +\widehat{m}(t) \leq y} \leq \Cr{cst:ub2} (z\vee1)\, {\rm e}^{\sqrt{2}y-\frac{z^2}{4\sqrt{t}\log t}} + o(1)
        \end{equation}
        for all $t$ large enough, where $o(1)$ denotes a function converging to $0$ as $t\to\infty$, independent of the choice of $z$.
        \item There exists a constant $\Cl[c]{cst:lb21} > 0$ such that if $z\in [0,t^{3/8}]$ and again $y=o(\log t)$, then
    	\begin{equation} \label{eq:etaLargeLB}
    		\Pbe{z}{\widehat{\xi}}{M_t + \widehat{m}(t) \leq y} \geq \Cr{cst:lb21} \bigl( 1 - \exp\bigl\{-\Cr{cst:lb21} (z\vee1) \,{\rm e}^{\sqrt{2}y-\frac{z^2}{2\sqrt{t}}}\bigr\} \bigr)
    	\end{equation}
    	for all $t$ large enough.
    \end{enumerate}
\end{theorem}
In the proofs, we will repeatedly use upper moderate deviation probabilities of standard BBM. In the F-KPP wording of \cite[Proposition~2.1]{BoHa-18} (see also \cite{LC25}) and in combination with \eqref{eq:mckMin}, we obtain that for any $x=x(t)>0$ fulfilling $\lim_{t\to\infty}x(t)=\infty$ and $x = o(t/\log t)$, we have
\begin{equation}\label{eq:BramsonAsymp}
	 	\P_{0}^{1}\Big(M_{t} \leq -\sqrt{2}t + \frac{3}{2\sqrt{2}} \log (t)- x \Big) \sim \Cl[c]{cst:Mt} x {\rm e}^{-\sqrt{2}x-\frac{x^2}{2t}}	\quad\text{ as } t\to\infty,
\end{equation}
for some universal constant $\Cr{cst:Mt}\in(0,\infty)$; here, we stick to the standard notation of writing $x(t) \sim y(t)$ as $t\to\infty$ if $y(t) \ne 0$ for all $t$ large enough and  $ \lim_{t\to\infty}\frac{x(t)}{y(t)} = 1$. Corollary~1.2 in \cite{LC25} also gives that the convergence in \eqref{eq:BramsonAsymp} is uniform for $x\in[a_t,b_t]$, for any functions $a_t\leq b_t$ with $a_t,b_t\xrightarrow{t\to\infty}\infty$ and $a_t,b_t=o(t/\log t)$.

The proofs of Theorems~\ref{the:min_case1} and \ref{the:min_case2} are organized as follows: In Section~\ref{sec:UB}, we derive an upper bound for probabilities of the form $\Pbe{z}{\widehat{\xi}}{M_t \leq y}$ in terms of moments of a certain sequence of stopping times and probabilities for homogeneous BBM. In Section~\ref{sec:LB} we modify this sequence to obtain a corresponding lower bound. Section~\ref{sec:tau_mom} and its subsections form the technical core of this article: we derive the required moment estimates for the number of entries of the aforementioned sequences in a unit interval. These are then combined with the asymptotics for moderate deviation probabilities of homogeneous BBM given in \eqref{eq:BramsonAsymp}. The latter provide the probabilities associated with the individual entries of the stopping-time sequences, while the moment estimates control how many such entries contribute.
In Section~\ref{sec:proof_bbm_est} we finally put these results together to prove Theorems~\ref{the:min_case1} and \ref{the:min_case2}.
In Section~\ref{sec:proof_bbm_main}, we prove Theorem~\ref{thm:main_bbm}. The proof uses an amplification argument to deduce tightness from Theorems~\ref{the:min_case1} and \ref{the:min_case2}.
Finally, in Section~\ref{sec:ub_below_zero}, we establish an upper bound for the probability $\Pbe{-z}{\widehat{\xi}}{M_t \leq y}$, $z \ge 0,$ i.e.\ when the process starts below the origin. This estimate is used later in the proofs of our main results concerning the $\varepsilon$-fronts.

\subsection{Statement of upper bound} \label{sec:UB}
We begin by rewriting the probability
\(  \Pbe{z}{\widehat{\xi}}{M_t \leq y}\).
For this purpose, we introduce the following suitable coloring scheme: the original particle of BBM$\widehat{\xi}$ is black. If an additional descendant of a black particle is generated as a result of the increased branching rate in the positive half-line, then that descendant is colored red; more precisely, while a black particle is at site $x\in\R$, any additional branching event, occurring at rate $(\widehat{\xi}(x)-1)$, produces an offspring configuration in which one particle remains black and all additional particles, if any, are colored red. In particular, this implies that black particles produce red particles at rate $(1-p_1)(\widehat{\xi}(x)-1)$.
All other descendants of black particles remain black. Moreover, all descendants of red particles are also colored red.
Whenever a red particle hits the origin, that particle is recolored black (in particular, note that this entails that red particles do not enter the negative half-line). Once recolored, the particle again follows the same coloring scheme as above, i.e.\ it may again produce red particles through the additional branching, which in turn also switch their color upon hitting the origin. 

\begin{figure}[htbp]
\centering
\includegraphics[width=\textwidth]{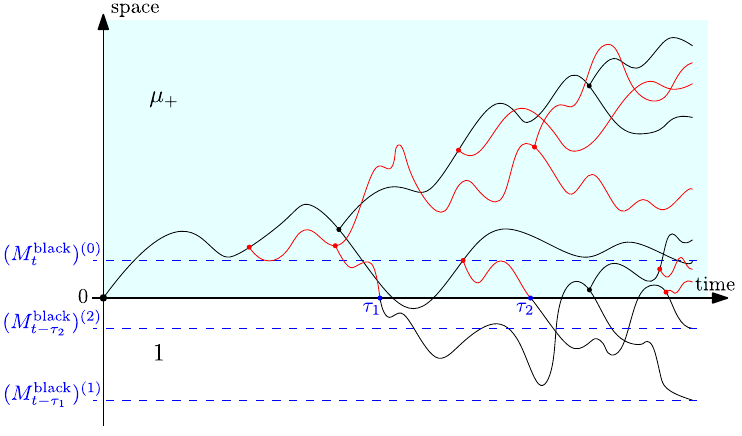}
\caption{Illustration of the colored process and its functionals.} 
\label{fig:tau}
\end{figure}

Next, we construct a sequence of stopping times together with a corresponding sequence of random variables describing the leftmost positions of certain sets of black particles.
Set $\tau_0=0$ and let $(M^{\text{black}}_t)^{(0)}$ denote the position of the leftmost black particle at time $t$ among the original black particle and its all-black descendants (i.e. \ those whose ancestry contains no intermediate red particle).
For $i\geq1$, let $\tau_i$ denote the $i$-th element of the (almost surely infinite) increasing sequence of first hitting times of $0$ by red particles. Equivalently, $\tau_i$ is the $i$-th time at which a red particle is recolored black. 
Let $(M^{\text{black}}_t)^{(i)}$ denote the position of the leftmost black particle, at time $\tau_i+t$, among the descendants of the particle recolored at time $\tau_i$ that have an all-black genealogy.
A visualization of this coloring scheme as well as of the sequences $\tau_i$ and $(M^{\text{black}}_t)^{(i)}$ is provided in Figure~\ref{fig:tau}.

Since black particles produce black descendants at constant rate $1$, each black particle, together with those of its descendants which have a completely black genealogy, evolves as a homogeneous BBM with branching rate $1$.
Consequently, we have
\begin{equation}\label{eq:black_min}
    \Pbe{z}{\widehat{\xi}}{(M^{\text{black}}_t)^{(i)} \leq y} = \begin{cases}
        \Pbe{z}{1}{M_t \leq y}, & i=0, \\
        \Pbe{0}{1}{M_t \leq y}, & i\geq1,
    \end{cases}
\end{equation}
for any $z,y\in\R$.

Note that the black particles \emph{fully} encode everything that happens \emph{below} the origin for the non-homogeneous BBM$\widehat{\xi}$.
Taking into account the time shifts $\tau_i$, we may represent the leftmost particle position of BBM$\widehat{\xi}$ as
\begin{equation}\label{eq:tau_construction}
    \Pbe{z}{\widehat{\xi}}{M_t \leq y} = \Pbe{z}{\widehat{\xi}}{\exists i\in\N_0 \, :\, {(M^{\text{black}}_{t-\tau_i})}^{(i)} \leq y}
\end{equation}
for any $z\in\R$ and $y<0$ (for the sake of notational simplicity, if $t\leq0$ we set $(M^{\text{black}}_t)^{(i)}=+\infty$ under $\mathbb{P}_z^{\widehat{\xi}}$).

This construction turns out to be a useful tool in order to further analyze the minimal position of BBM$\widehat \xi$, since the BBMs occurring on the right-hand side of \eqref{eq:tau_construction} are all homogeneous in nature; therefore, once we have an understanding of the $\tau_i$, we can take advantage of the profound understanding available for homogeneous BBM.

We define 
\begin{equation} \label{eq:Z}
Z_t(y) := \sum_{k=0}^{\floor{t}} \#\{i\geq1\,:\,\tau_i \in [k,k+1)\}\, \Pbe{0}{1}{M_{t-k} \leq y}
\end{equation}
and obtain the following simple lemma.
\begin{lemma}\label{lem:prob_ub}
	Let $y< 0$ and $z\in\R$, then
	\begin{equation*}
		\Pbe{z}{\widehat{\xi}}{M_t \leq y} \leq \Pbe{z}{1}{M_{t} \leq y} + \Eeb{z}{\widehat{\xi}}{Z_t(y)}.
	\end{equation*}
\end{lemma}
\begin{proof}
	Using \eqref{eq:tau_construction} and a union bound,
    \begin{align*}
        \Pbe{z}{\widehat{\xi}}{M_t \leq y} &= \Pbe{z}{\widehat{\xi}}{\exists i\in\N_0 \, :\, {(M^{\text{black}}_{t-\tau_i})}^{(i)} \leq y} \\
        &\leq \Pbe{z}{\widehat{\xi}}{ {(M^{\text{black}}_{t})}^{(0)} \leq y} + \sum_{i\geq1} {\Pbe{z}{\widehat{\xi}}{ {(M^{\text{black}}_{t-\tau_i})}^{(i)} \leq y}}.
    \end{align*}
    Then, conditioning on $\tau_i$ while discretizing time, the second summand is upper bounded by
    \begin{align*}
        \sum_{i\geq 1} \sum_{k=0}^{\floor{t}} \Ee{z}{\widehat{\xi}}{ \1{\tau_i \in [k,k+1)} \Pbe{z}{\widehat{\xi}}{ {(M^{\text{black}}_{t-\tau_i})}^{(i)} \leq y \mid \tau_i } }.
    \end{align*}
    The strong Markov property together with \eqref{eq:black_min} yields
    \begin{align*}
        \Pbe{z}{\widehat{\xi}}{M_t \leq y} &\leq \Pbe{z}{1}{M_{t} \leq y} + \sum_{k=0}^{\floor{t}} \Ee{z}{\widehat{\xi}}{ \sum_{i\geq 1} \1{\tau_i \in [k,k+1)} } \, \Pbe{0}{1}{ {M}_{t-k} \leq y }\\
        &= \Pbe{z}{1}{M_{t} \leq y} + \Ee{z}{\widehat{\xi}}{Z_t(y)},
    \end{align*}
    where to get the inequality we used the monotonicity of $\Pbe{0}{1}{ {M}_{t} \leq y }$ in $t$. The claim follows.
\end{proof}

\subsection{Statement of lower bound} \label{sec:LB}

A difficulty with the sequence $\tau = (\tau_i)_{i \in \N}$ is that the paths of the black particles are not independent of $\tau$. Indeed, while above the origin, black particles may produce red particles, which in turn generate new stopping times $\tau_i$; this dependence prevents the straightforward application of a second moment method for obtaining a lower bound. 
To get around this problem, we will disregard the branching of some black particles by not letting them produce red particles. More precisely, we do so by ``switching off'' the increased branching rate for them.

Recall the process of colored particles from Section \ref{sec:UB}.
In order to accommodate the above-mentioned switching off, similarly to the construction for the upper bound, we now define a modified sequence of stopping times $\tau^\prime = (\tau_i^\prime)$.
As before, the stopping times are given by the successive first hitting times of the origin by red particles. The only difference compared to the previous construction is that those red particles that hit $0$ are recolored {\em blue} instead of black. Blue particles evolve as homogeneous BBM with branching rate $1$, and their offspring is blue as well (in particular, they do not generate any red particles anymore).
Hence, $\tau'_i$ is the $i$-th time at which a red particle is recolored blue.

Analogously to Section~\ref{sec:UB}, for $i\geq1$ we let $(M^{\text{blue}}_t)^{(i)}$ denote the leftmost blue particle position (among the particle recolored at time $\tau'_i$ and all of its blue descendants) when time $t$ has elapsed since time  $\tau'_i$.

\begin{figure}[htbp]
\centering
\includegraphics[width=\textwidth]{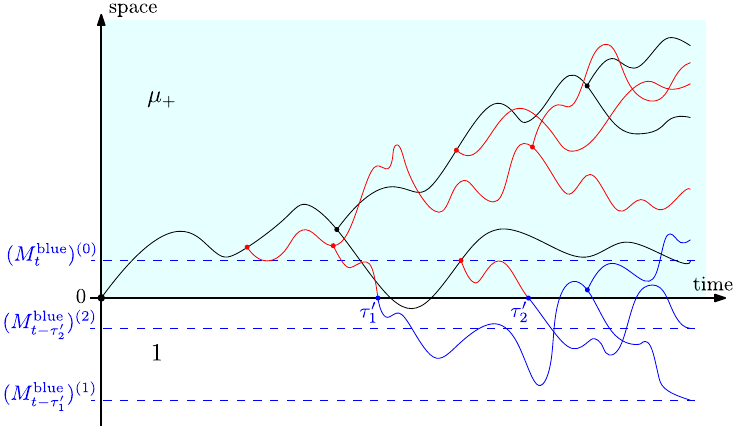}
\caption{A sketch of the sequence $\tau_i'$. As in Figure \ref{fig:tau}, but with blue particles that do not generate any red particles.} 
\label{fig:tauPrime}
\end{figure}

This modified construction removes the dependence described above and the consideration of such particles and the respective $\tau'$ will facilitate the proof of Lemma \ref{lem:prob_lb}; see Figure \ref{fig:tauPrime} for an illustration.

Under the coupling of the same branching and evolution dynamics, every entry of $\tau^\prime$ is also contained in $\tau$. Consequently, for every $z\in\R$ and $y<0$, one has
\begin{equation}\label{eq:tauPrime_construction}
    \Pbe{z}{\widehat{\xi}}{M_t \leq y} \geq \Pbe{z}{\widehat \xi}{\exists i\in\N \, :\, (M^{\text{blue}}_{t-\tau^\prime_i})^{(i)} \leq y},
\end{equation}
similarly to \eqref{eq:tau_construction}.
Furthermore, recalling $r^*$ from \eqref{def:r*} we set 
\begin{equation} \label{eq:IDef}
\begin{split}
    \Lambda^\prime\coloneqq& \big[\lfloor\sqrt{t}\vee (r^*t -\sqrt{t})\rfloor, \lceil(2\sqrt{t})\vee(r^*t +\sqrt{t})\rceil\big], \\ \Lambda \coloneqq& \big[\lfloor\sqrt{t}\vee (r^*t -\sqrt{t})\rfloor, \lceil(2\sqrt{t})\vee(r^*t +\sqrt{t})\rceil-1\big] \cap\N;
\end{split}
\end{equation}
in this context, note that for $\eta_0\in[1/2,\infty]$ one has $r^*t<\sqrt{t}$, while for $\eta_0\in(0,1/2)$ we infer $\frac{r^*t}{\sqrt{t}} \to \infty$ as $t \to \infty$.
Now for $t$ large enough so that $\Lambda^\prime \subset [0,t]$, define 
\begin{equation} \label{eq:Z'}
Z^\prime_t(y) \coloneqq \sum_{k\in \Lambda}\# \{i : \tau^\prime_i \in [k,k+1)\} \, \Pbe{0}{1}{M_{t-(k+1)} \leq y }.
\end{equation}
We get the following lower bound.
\begin{lemma}\label{lem:prob_lb}
    Let $y<0$ and $z\geq 0$. Then, for $t$ large enough, 
    \begin{align*}
        \Pbe{z}{\widehat{\xi}}{M_t \leq y} \geq  \mathbb{P}_{z}^{\widehat{\xi}}\Bigl(Z^\prime_t(y)\geq \frac12 \mathbb{E}_{z}^{\widehat{\xi}}[Z^\prime_t(y)]\Bigr) \bigl( 1 - {\rm e}^{-\frac{1}{2} \mathbb{E}_{z}^{\widehat{\xi}}[Z^\prime_t(y)]} \bigr).
    \end{align*}
\end{lemma}
\begin{proof}
    By \eqref{eq:tauPrime_construction} we have
    \begin{align*}
        \Pbe{z}{\widehat{\xi}}{M_t \leq y} &\geq \Pbe{z}{\widehat{\xi}}{\exists \tau^\prime_i\in \Lambda^\prime \, :\, (M^{\text{blue}}_{t-\tau^\prime_i})^{(i)} \leq y} = 1 - \Ee{z}{\widehat \xi}{\prod_{\substack{i\in\N \\ \tau^\prime_i\in \Lambda^\prime}}\Pbe{z}{\widehat \xi}{(M^{\text{blue}}_{t-\tau^\prime_i})^{(i)} > y\, \big\vert\, \tau^\prime}}  \\
        &= 1 - \Ee{z}{\widehat \xi}{\prod_{k\in \Lambda}\prod_{\tau^\prime_i \in [k,k+1)}\Pbe{z}{\widehat \xi}{(M^{\text{blue}}_{t-\tau^\prime_i})^{(i)}  > y \,\big\vert\, \tau^\prime_i}}.
    \end{align*}
    Here, since we only track blue particles that do not generate any red offspring, and hence cannot create a new entry in the sequence $(\tau^\prime_i)_i$, each $(M^{\text{blue}}_t)^{(i)}$ is independent of any of the future stopping times $\tau^\prime_j$ for $j>i$.
    Then, given $\tau_i^\prime \in [k,k+1)$, we estimate 
    \begin{align*}
        \Pbe{z}{\widehat \xi}{(M^{\text{blue}}_{t-\tau^\prime_i})^{(i)}  > y \,\big\vert\, \tau^\prime_i} \leq \Pbe{0}{1}{M_{t-(k+1)} > y }.
    \end{align*}
     Plugging this into the previous display yields
    \begin{equation}\label{eq:Mt_lb1}
    \begin{split}
        \Pbe{z}{\widehat{\xi}}{M_t \leq y} &\geq 1 - \Ee{z}{\widehat{\xi}}{\prod_{k\in \Lambda}\Bigl( 1 - \Pbe{0}{1}{M_{t-(k+1)} \leq y } \Bigr)^{\# \{i: \tau_i^\prime \in [k,k+1)\} } } \\
        &\geq 1 - \Ee{z}{\widehat{\xi}}{ \exp\Big\{- \sum_{k\in \Lambda}\# \{i : \tau_i^\prime \in [k,k+1)\} \, \Pbe{0}{1}{M_{t-(k+1)} \leq y }\Big\}}\\
        &= 1 - \Eeb{z}{\widehat{\xi}}{ \exp\{-Z^\prime_t(y)\}},
    \end{split}
    \end{equation}
    where we took advantage of the inequality $1-x \leq {\rm e}^{-x}$ for all $x\in\R$ in the second inequality. For any non-negative random variable $X$ we can estimate
    \begin{align*}
        \mathbb{E}\bigl[ {\rm e}^{-X} \bigr] &= \mathbb{E}\bigl[ {\rm e}^{-X} \1{X\geq\frac{1}{2}\mathbb{E}[X]} \bigr] + \mathbb{E}\bigl[ {\rm e}^{-X} \1{X< \frac{1}{2}\mathbb{E}[X]} \bigr] \\
        &\leq {\rm e}^{-\frac{1}{2}\mathbb{E}[X]}\mathbb{P}\bigl(X\geq\mathbb{E}[X]/2\bigr) + \bigl( 1 - \mathbb{P}\bigl(X\geq\mathbb{E}[X]/2\bigr) \bigr) \\
        &= 1 - \mathbb{P}\bigl(X\geq\mathbb{E}[X]/2\bigr) \bigl( 1 - {\rm e}^{-\frac{1}{2}\mathbb{E}[X]} \bigr);
    \end{align*}
    together with \eqref{eq:Mt_lb1} this finishes the proof.
\end{proof}
The next lemma supplies a sufficient condition for the probability $\mathbb{P}_{z}^{\widehat{\xi}}\bigl(Z^\prime_t(y)\geq\frac{1}{2} \mathbb{E}_{z}^{\widehat{\xi}}[Z^\prime_t(y)]\bigr)$ in Lemma~\ref{lem:prob_lb} being bounded away from zero.
\begin{lemma}\label{lem:paley}
	In the setting of Lemma~\ref{lem:prob_lb}, assume there exists a constant $C>0$ such that
	\begin{align}\label{cond:paley}
		\mathbb{E}_{z}^{\widehat{\xi}}\Bigl[\bigl(\# \{i:\tau_i^\prime \in [k,k+1)\}\bigr)^2\Bigr] \leq C\mathbb{E}_{z}^{\widehat{\xi}}\Bigl[\# \{i:\tau_i^\prime \in [k,k+1)\}\Bigr]^2
	\end{align}
    for all $k \in \Lambda$.
	Then
	\begin{align*}
		\mathbb{P}_{z}^{\widehat{\xi}}\Bigl(Z^\prime_t(y)\geq\frac{1}{2} \mathbb{E}_{z}^{\widehat{\xi}}[Z^\prime_t(y)]\Bigr) \geq \frac{1}{4C}.
	\end{align*}
\end{lemma}
\begin{proof}
	We use the Paley--Zygmund inequality (see e.g.\ Theorem 2.3.1 in \cite{MR5088005}) to obtain
    \begin{align*}
\mathbb{P}_{z}^{\widehat{\xi}}\Bigl(Z^\prime_t(y)\geq\frac{1}{2} \mathbb{E}_{z}^{\widehat{\xi}}[Z^\prime_t(y)]\Bigr) \geq \frac{1}{4} \frac{\bigl(\mathbb{E}_{z}^{\widehat{\xi}}[Z^\prime_t(y)]\bigr)^2}{\mathbb{E}_{z}^{\widehat{\xi}}[Z^\prime_t(y)^2]}
    \end{align*}
    and with Cauchy--Schwarz it follows that
    \begin{align*}
        \mathbb{E}_{z}^{\widehat{\xi}}\bigl[\#\{i:\tau_i^\prime \in [k,k+1)\} &\#\{i:\tau_i^\prime \in [l,l+1)\}\bigr]\\ 
        &\leq \Bigl(\mathbb{E}_{z}^{\widehat{\xi}}\bigl[(\#\{i:\tau_i^\prime \in [k,k+1)\})^2\bigr] \mathbb{E}_{z}^{\widehat{\xi}}\bigl[(\#\{i:\tau_i^\prime \in [l,l+1)\})^2\bigr]\Bigr)^{1/2}
    \end{align*}
    from which we infer
    \begin{align*}
\mathbb{E}_{z}^{\widehat{\xi}}[Z^\prime_t(y)^2]
        &= \sum_{k,l\in \Lambda}\mathbb{E}_{z}^{\widehat{\xi}}\bigl[\#\{i:\tau_i^\prime \in [k,k+1)\} \#\{i:\tau_i^\prime \in [l,l+1)\}\bigr]\\ &\qquad\qquad\quad\cdot \Pbe{0}{1}{M_{t-(k+1)} \leq y } \Pbe{0}{1}{M_{t-(l+1)} \leq y } \\
        &\leq \Bigl\{ \sum_{k\in \Lambda} \Bigl(\mathbb{E}_{z}^{\widehat{\xi}}\Bigl[\bigl(\#\{i:\tau_i^\prime \in [k,k+1)\}\bigr)^2\Bigr]\Bigr)^{1/2} \, \Pbe{0}{1}{M_{t-(k+1)} \leq y } \Bigr\}^2.
    \end{align*}
    Plugging the upper bound of  condition \eqref{cond:paley} into the right-hand side yields
$\mathbb{E}_{z}^{\widehat{\xi}}[Z^\prime_t(y)^2] \leq C \bigl(\mathbb{E}_{z}^{\widehat{\xi}}[Z^\prime_t(y)]\bigr)^2,
    $
    which finishes the proof.
\end{proof}

Lemmas~\ref{lem:prob_ub} to \ref{lem:paley} suggest that in order to further bound $\Pbe{z}{\widehat{\xi}}{M_t \leq y}$ from above and below, we need to upper bound the first and second moments, as well as lower bound the first moments, of the number of stopping times $\tau_i$ and $\tau^\prime_i$ taking values in an interval $[k,k+1)$. We will do so in the subsequent sections.

\subsection{Upper and lower bounds  of the moments of the $\tau_i$s in intervals}\label{sec:tau_mom}
In this subsection we derive the results needed for applying the Paley--Zygmund inequality in an effective way.

\subsubsection{Upper bounding the expected number of $\tau_i$s in intervals}
We define for $t \ge 0$ the first hitting time 
\begin{equation*}
    H_0^{\geqslant t} \coloneqq \inf\{s\geq t \, : \, B_s=0 \}
\end{equation*}
    of $0$ by Brownian motion after time $t$, and abbreviate $H_0\coloneqq H_0^{\geqslant 0}$.
We start with an auxiliary result that will be used multiple times throughout this section.
\begin{lemma}\label{lem:num_tau_aux}
    There exists a constant $C\in(0,\infty)$ such that for $k\geq0$ and $z\geq0$,
    \begin{equation}\label{eq:num_tau_aux-b1}
        \E{z}{{\rm e}^{\int_0^{k+1} \widehat{\xi}(B_r)\, \mathrm{d}r} \1{\exists s\in [k,k+1) \, : \, B_s = 0} } \leq C(z\vee 1)\,{\rm e}^{-\frac{z^2}{2(k+1)}}(k\vee 1)^{-3/2}  {\rm e}^{\mu_+ k}
    \end{equation}
    and
    \begin{equation}\label{eq:num_tau_aux-b2}
        \E{-z}{{\rm e}^{\int_0^{k+1} \widehat{\xi}(B_r)\, \mathrm{d}r} \1{\exists s\in [k,k+1) \, : \, B_s = 0} } \leq C\,{\rm e}^{-z}(k\vee 1)^{-3/2}  {\rm e}^{\mu_+ k}.
    \end{equation}
\end{lemma}
\begin{proof}
    We first prove \eqref{eq:num_tau_aux-b1}. 
    For $k=0$ and since $z \ge 0$, by the reflection principle for Brownian motion and standard Gaussian estimates (\cite[Chapter 1.2]{AdTa-07}), we can simply upper bound this expectation by ${\rm e}^{\mu_+} 2\Pb{z}{B_1 \leq 0} \leq C{\rm e}^{\mu_+} {\rm e}^{-z^2/2}$. For $k \ge 1$, we can upper bound 
	\begin{align*}
		&\E{z}{{\rm e}^{\int_0^{k+1} \widehat{\xi}(B_r)\, \mathrm{d}r} \1{\exists s\in [k,k+1) \, : \, B_s = 0} } \leq {\rm e}^{\mu_+} \E{z}{{\rm e}^{\int_0^{H_0^{\geqslant k}} \widehat{\xi}(B_r)\, \mathrm{d}r} \1{H_0^{\geqslant k}\in [k,k+1)} }.
	\end{align*}	
	Disintegrating according to the value of the first hitting time of $0$ and  writing $\Pb{z}{H_0\in\mathrm{d}s}$ for the density at $s \in (0,\infty)$ of $\P_z \circ H_0^{-1}$ with respect to Lebesgue measure,  the strong Markov property of Brownian motion supplies us with
	\begin{align*}
		&\E{z}{{\rm e}^{\int_0^{H_0^{\geqslant k}} \widehat{\xi}(B_r)\, \mathrm{d}r} \1{H_0^{\geqslant k}\in [k,k+1)}} \\
        &\quad = \int_0^k {\rm e}^{\mu_+ s}\,\Pb{z}{H_0\in\mathrm{d}s}\,\E{0}{{\rm e}^{\int_0^{H_0^{\geqslant k-s}} \widehat{\xi}(B_r)\, \mathrm{d}r} \1{H_0^{\geqslant k-s}\in [k-s,k-s+1)}} \, \mathrm{d}s\\
        &\phantom{\quad =}+ \int_k^{k+1} {\rm e}^{\mu_+ s} \,\Pb{z}{H_0\in\mathrm{d}s} \, \mathrm{d}s.
	\end{align*}
    The event $\{H_0^{\geqslant k-s}\in [k-s,k-s+1)\}$ implies that the last passage time of the origin during the time interval $[0,k-s+1)$ also happens during $[k-s,k-s+1)$. Upper bounding the functional of $B_r$ up to the last passage time, we get a functional pertaining to a standard Brownian bridge (cf. \cite[Part I IV.19]{handbookBM}). Letting $\zeta^T$ denote a standard Brownian bridge of duration $T$ and taking advantage of the fact that for $u \in [k-s,k+1-s]$ the expression $\E{0}{{\rm e}^{\int_0^{u} \widehat{\xi}(\zeta^{u}(r)) \, \mathrm{d}r}}$ is maximized for $u = k+1-s$, we obtain the upper bound
    \begin{align*}
		&\E{0}{{\rm e}^{\int_0^{H_0^{\geqslant k-s}} \widehat{\xi}(B_r)\, \mathrm{d}r} \1{H_0^{\geqslant k-s}\in [k-s,k+1-s)}} \\
        &\quad\leq \E{0}{{\rm e}^{\int_0^{k+1-s} \widehat{\xi}(\zeta^{k+1-s}(r)) \, \mathrm{d}r}} \Pb{0}{H_0^{\geqslant k-s}\in [k-s,k+1-s)},
	\end{align*}
    for $s\in[0,k]$.
    Since the time a standard Brownian bridge spends above zero is uniformly distributed (see \cite{Pi-99}), we obtain
    \begin{align*}
        \E{0}{{\rm e}^{\int_0^{k+1-s} \widehat{\xi}(\zeta^{k+1-s}(r)) \, \mathrm{d}r}} &\leq {\rm e}^{k+1-s} \frac{{\rm e}^{(\mu_+ -1) (k+1-s)}}{(\mu_+ -1)(k+1-s)}=\frac{{\rm e}^{\mu_+ (k+1-s)}}{(\mu_+ -1)(k+1-s)}
    \end{align*}
    and therefore, combining the three displays above, 
    \begin{equation}\label{eq:num_tau_ub_b1}
    \begin{split}
    		&\E{z}{{\rm e}^{\int_0^{H^{\geqslant k}_0} \widehat{\xi}(B_r)\, \mathrm{d}r} \1{H_0^{\geqslant k}\in [k,k+1)}} \\ 
            &\leq C {\rm e}^{\mu_+ k} \int_0^k \frac{\Pb{0}{H_0^{\geqslant k-s}\in [k-s,k+1-s)}}{k+1-s}  \Pb{z}{H_0\in\mathrm{d}s} \, \mathrm{d}s \\
    	&\phantom{\leq} + {\rm e}^{\mu_+ (k+1)}\, \Pb{z}{H_0\in[k,k+1)}.
    \end{split}
    \end{equation}
    To further bound the right-hand side of \eqref{eq:num_tau_ub_b1} from above, we first estimate $\Pb{0}{H_0^{\geqslant l}\in [l,l+1)}$ for any $l>0$.
    Observe that due to the Markov property, the reflection principle implies
    \begin{align*}
        &\Pb{0}{H_0^{\geqslant l}\in [l,l+1)} = 2 \E{0}{\Pb{y}{B_{1} \leq 0}_{y= \abs{B_l}}} \\
        &= 2\E{0}{\Pb{y}{B_{1} \leq 0}_{y= \abs{B_l}}\1{\abs{B_l}>\sqrt{l}}} + 2\E{0}{\Pb{y}{B_{1} \leq 0}_{y= \abs{B_l}}\1{\abs{B_l}\leq \sqrt{l}}}.
    \end{align*}
    The first summand can be upper bounded by
    \begin{align*}
        2\Pb{0}{B_1 \geq \sqrt{l}} \Pb{0}{|B_l|>\sqrt{l}} \leq \frac{C}{\sqrt{l}},
    \end{align*}
        and the second summand by
    \begin{align*}
        2\int_{-\sqrt{l}}^{+\sqrt{l}} \Pb{0}{B_1 \geq \abs{x}} \frac{1}{\sqrt{2\pi l}} {\rm e}^{-\frac{x^2}{2l}} \, \mathrm{d}x \leq \frac{C}{\sqrt{ l}};
    \end{align*}
     combining these bounds provides the total upper bound 
    \begin{align*}
    		\Pb{0}{H_0^{\geqslant l}\in [l,l+1)} \leq \frac{C}{\sqrt{l}}.
    \end{align*}
    Putting the above together, we can upper bound the integral in \eqref{eq:num_tau_ub_b1} by
    \begin{align}\label{eq:expUBIntMaxProb}
    \begin{split}
    		C \int_0^k \frac{\Pb{z}{H_0\in\mathrm{d}s}}{(k+1-s)\sqrt{k-s}} \, \mathrm{d}s \leq \frac{C}{k^{3/2}} \Pb{z}{H_0\leq k/2} + \frac{Cz}{k^{3/2}} {\rm e}^{-\frac{z^2}{2k}} \int_{k/2}^k \frac{1}{(k+1-s)\sqrt{k-s}} \, \mathrm{d}s
    \end{split}
    \end{align}
         where we split the integral into two parts and used the explicit density of $H_0$ (provided in Lemma~\ref{lem:HittingTimeDensity} for the sake of completeness).
    
    The remaining integral in the second summand of \eqref{eq:expUBIntMaxProb} is of constant order. To further upper bound the first summand on the right-hand side, we again use the reflection principle  as well as the scale invariance and reversibility of Brownian motion. We infer that $\Pb{z}{H_0\leq k/2} = 2\Pb{0}{B_1 \geq z/\sqrt{k/2}} \leq {\rm e}^{-z^2/k}$ via the simple estimate
    \begin{equation*}
        \int_{0}^\infty \frac{1}{\sqrt{2\pi}} {\rm e}^{-\frac{(x+z)^2}{2}} \,\mathrm{d}x = {\rm e}^{-z^2/2} \int_{0}^\infty \frac{1}{\sqrt{2\pi}} {\rm e}^{-x^2/2 - xz} \,\mathrm{d}x \leq \frac{1}{2} {\rm e}^{-z^2/2}
    \end{equation*}
    for $z\geq0$. Altogether, this provides an upper bound for the first summand in \eqref{eq:num_tau_ub_b1}.
    
    The probability in the second summand of \eqref{eq:num_tau_ub_b1} can again be upper bounded by appealing to  Lemma \ref{lem:HittingTimeDensity}; combining the resulting estimate with a separate consideration of the case $k=0$, Equation~\eqref{eq:num_tau_aux-b1} follows.

    Next, we prove \eqref{eq:num_tau_aux-b2}.
    Conditioning on the first hitting time of $0$ together with the strong Markov property,
    \begin{align*}
        &\E{-z}{{\rm e}^{\int_0^{k+1} \widehat{\xi}(B_r)\, \mathrm{d}r} \1{\exists s\in [k,k+1) \, : \, B_s = 0} }\\ 
        &= \E{-z}{{\rm e}^{H_0}\1{H_0\leq k+1}\E{0}{{\rm e}^{\int_0^{k+1-h} \widehat{\xi}(B_r)\, \mathrm{d}r} \1{\exists s\in [0\vee (k-h),k-h+1) \, : \, B_s = 0} }_{h=H_0}}.
    \end{align*}
    By \eqref{eq:num_tau_aux-b1} we may upper bound the inner expectation by
    \begin{equation*}
        \E{0}{{\rm e}^{\int_0^{k+1-h} \widehat{\xi}(B_r)\, \mathrm{d}r} \1{\exists s\in [k-h,k-h+1) \, : \, B_s = 0} } \leq C ((k-h)\vee1)^{-3/2} {\rm e}^{\mu_+(k-h)}.
    \end{equation*}
    Hence,
    \begin{equation}\label{eq:ub_b0-eq2}
    \begin{split}
        &\E{-z}{{\rm e}^{\int_0^{k+1} \widehat{\xi}(B_r)\, \mathrm{d}r} \1{\exists s\in [k,k+1) \, : \, B_s = 0} } \leq C {\rm e}^{\mu_+k}\,\E{-z}{((k-H_0)\vee1)^{-\frac32} {\rm e}^{-(\mu_+-1)H_0}}.
    \end{split}
    \end{equation}
    For some $a>0$ we now split the expectation into two parts, where $H_0\leq a$ or $H_0>a$. For the case $H_0\leq a$ we upper bound
    \begin{align*}
        &\E{-z}{\1{H_0 \leq a}((k-H_0)\vee1)^{-\frac32} {\rm e}^{-(\mu_+-1)H_0}} \leq ((k-a)\vee1)^{-3/2}\,\E{-z}{{\rm e}^{-H_0/2}}.
    \end{align*}
    In the case $H_0>a$ we upper bound
    \begin{align*}
        &\E{-z}{\1{H_0 > a}((k-H_0)\vee1)^{-\frac32} {\rm e}^{-(\mu_+-1)H_0}} \leq {\rm e}^{-a/2}\E{-z}{{\rm e}^{-H_0/2}}.
    \end{align*}
    For the expectation, \cite[p. 205, 2.0.1]{handbookBM} gives
    \[
        \E{-z}{{\rm e}^{-H_0/2}} = {\rm e}^{-z}.
    \]
    Choosing $a=k/2$ together with \eqref{eq:ub_b0-eq2} yields \eqref{eq:num_tau_aux-b2}.
\end{proof}
\begin{lemma}\label{lem:num_tau_ub}
	There exists a constant $C \in (0,\infty)$ such that for $z \geq 0$ and $k\geq0$:
	\begin{align*}
		\Ee{z}{\widehat{\xi}}{\#\{i\geq1\,:\,\tau_i \in [k,k+1)\}} \leq C(z\vee 1){\rm e}^{-\frac{z^2}{2(k+1)}}(k\vee 1)^{-3/2}  {\rm e}^{\mu_+ k}.
	\end{align*}
\end{lemma}
\begin{proof}
	Recalling the notation from Section \ref{sec:conn}, we start with the coarse upper bound
	\begin{align}\label{eq:num_tau_ub_b0}
		\Ee{z}{\widehat{\xi}}{\#\{i\geq1\,:\,\tau_i \in [k,k+1)\}} \leq \Ee{z}{\widehat{\xi}}{\sum_{\nu\in N(k+1)} \1{\exists s\in [k,k+1) \, : \, X^\nu_s = 0} };
	\end{align}	
indeed, this inequality holds true since each $\tau_i \in [k,k+1)$  can be mapped injectively to (at least) one black particle in $N(k+1)$ (or even to its entire black progeny in  $N(k+1)$) since its black genealogy is non-trivial at any positive time due to \eqref{eq:mom_offspr} which ensures that at each branching event at least one (black) particle is produced.
By the many-to-one lemma (see Lemma~\ref{lem:many-to-one}), the right-hand side of \eqref{eq:num_tau_ub_b0} is equal to
	\begin{align} \label{eq:FKeq}
		\E{z}{{\rm e}^{\int_0^{k+1} \widehat{\xi}(B_r)\, \mathrm{d}r} \1{\exists s\in [k,k+1) \, : \, B_s = 0} }.
	\end{align}
Applying Lemma~\ref{lem:num_tau_aux} then finishes the proof.
\end{proof}

\subsubsection{Lower bounding the expected number of $\tau_i'$s in intervals}
The next lemma establishes a lower bound for the expectation of the number of $\tau_i'$s in a unit interval. For this purpose, let
\begin{equation}\label{eq:runMin}
    \underline{X}_t \coloneqq \inf_{s\in[0,t]} X_s
\end{equation}
denote the running infimum of a process $(X_t)$. 
\begin{lemma}\label{lem:num_tau_lb}
	There exists a constant $C > 0$ such that for all $z\geq 0$ and $k\geq2$,
	\begin{align*}
		\Ee{z}{\widehat{\xi}}{\#\{i:\tau^\prime_i \in [k,k+1)\}} \geq C(z\vee 1){\rm e}^{-\frac{z^2}{2k}} k^{-3/2} {\rm e}^{\mu_+ k}.
	\end{align*}
\end{lemma}

\begin{proof}
    We start with observing that 
    \begin{equation*}
     \#\{i:\tau^\prime_i \in [k,k+1)\} \geq \sum_{\nu\in N(k)} \1{\nu \text{ is red}} \, \1{\text{descendant of }\nu\text{ hits }0\text{ in }[k,k+1)}.
    \end{equation*}
    We further lower bound the right-hand side by considering only the descendants of the first red particle. 
    Let $\sigma$ be the time of the first branching event of the entire BBM process.
    Moreover, let $\mathcal{A}$ denote the event that the first branching event leads to at least one red particle and that $\sigma\in[1,2]$. This event has positive probability uniform for $z\geq0$. Using the strong Markov property at time $\sigma$, and since any descendant of the first red particle that remained above the origin until time $k$ is itself red, we obtain
    \begin{equation}\label{eq:num_tauPrime_lb_b2}
    \begin{split}
        &\Ee{z}{\widehat{\xi}}{\#\{i:\tau^\prime_i \in [k,k+1)\}} \\
        &\geq \Ee{z}{\widehat{\xi}}{\mathds{1}_{\mathcal{A}}\Ee{B_{\sigma}}{\widehat{\xi}}{\sum_{\nu\in N(k-s)}  \1{\underline{X}^\nu_{k-s}>0}\1{\text{descendant of }\nu\text{ hits }0\text{ in }[k-s,k+1-s)}}_{s=\sigma}},
    \end{split}
    \end{equation}
    where $B_\sigma$ denotes the position of the (unique) particle at the first branching event.
    Applying the many-to-one lemma (Lemma~\ref{lem:many-to-one}) and discarding any additional branching within the time interval $[k-s,k+1-s)$, we lower bound \eqref{eq:num_tauPrime_lb_b2} by
    \begin{equation}\label{eq:num_tauPrime_lb_b3}
    \begin{split}
        &{\rm e}^{\mu_+(k-2)}\Ee{z}{\widehat{\xi}}{\mathds{1}_{\mathcal{A}}{\mathbb{P}_{B_{\sigma}}\bigl({\underline{B}_{k-s}>0, \, \exists\, l\in[k-s,k+1-s)\,: \,B_l=0}\bigr)_{s=\sigma}}} \\
        &\geq C\,{\rm e}^{\mu_+k} \,\Pb{z}{B_r> 0 \,\forall r\in[1,k), \exists\, l\in[k,k+1) : B_l=0},
    \end{split}
    \end{equation}
    where to obtain the inequality we furthermore restricted the trajectory of the initial particle to remain positive on the time interval $[1,k)$.
    Finally, we observe that using Lemma~\ref{lem:HittingTimeDensity} we can infer the lower bound 
    \begin{equation*}
        \mathbb{P}_{z}\bigl({B_r>0 \,\forall\,r\in[1,k), \, \exists\, l\in[k,k+1)\,: \,B_l=0}\bigr) \geq \frac{C(z\vee1){\rm e}^{-\frac{z^2}{2k}}}{k^{3/2}}.
    \end{equation*}
      Plugging this last estimate into \eqref{eq:num_tauPrime_lb_b3}, 
    we obtain the desired result.
\end{proof}
Another key ingredient in verifying the assumptions of  Lemma~\ref{lem:paley} is the following upper bound on the second moment.
\begin{lemma}\label{lem:num_tau_2nd}
	There exists a constant $C > 0$ such that for all $z\geq0$ and $k\geq2$,
	\begin{equation*}
		\Ee{z}{\widehat{\xi}}{\bigl(\#\{i:\tau_i^\prime \in [k,k+1)\}\bigr)^2} \leq \frac{C(z\vee 1)^2 {\rm e}^{2\mu_+ k}}{k^3} \max\big\{{\rm e}^{-\frac{z^2}{k+\sqrt{k}}}, k^3 {\rm e}^{-\mu_+ \sqrt{k}}\big\}.
	\end{equation*}
\end{lemma}

\begin{proof}
	We can apply an argument analogous to the one employed in deriving \eqref{eq:num_tau_ub_b0} in order to upper bound 
    \begin{equation}\label{eq:num_tau_2nd_moments}
    \begin{split}
        &\mathbb{E}_{z}^{\widehat{\xi}}\Bigl[\bigl(\#\{i:\tau_i^\prime \in [k,k+1)\}\bigr)^2\Bigr]\\ 
        &\leq \Ee{z}{\widehat{\xi}}{\#\{i:\tau_i^\prime \in [k,k+1)\}} + \Ee{z}{\widehat{\xi}}{\sum_{\substack{\nu_1, \nu_2\in N(k+1) \\ \nu_1 \neq \nu_2}} \1{\exists\,l_1,l_2\in[k,k+1)\,:\,X^{\nu_1}_{l_1}=X^{\nu_2}_{l_2}=0}}.
    \end{split}
    \end{equation}
    Since the sequence $(\tau^\prime_i)$ forms a subsequence of the sequence $(\tau_i)$, we may upper bound the first summand via Lemma~\ref{lem:num_tau_ub} to obtain
    \begin{align}\label{eq:num_tau_2nd_firstmom}
        \Ee{z}{\widehat{\xi}}{\#\{i:\tau_i^\prime \in [k,k+1)\}} \leq \frac{C(z\vee 1){\rm e}^{-\frac{z^2}{2(k+1)}}{\rm e}^{\mu_+ k}}{k^{3/2}}.
    \end{align}
    By the many-to-two lemma (Lemma~\ref{lem:many-to-two}) and the Markov property we can upper bound the second summand as
\begin{equation}\label{eq:num_tau_2nd_b0}
    \begin{split}
        &\Ee{z}{\widehat{\xi}}{\sum_{\substack{\nu_1, \nu_2\in N(k+1) \\ \nu_1 \neq \nu_2}} \1{\exists\,l_1,l_2\in[k,k+1)\,:\,X^{\nu_1}_{l_1}=X^{\nu_2}_{l_2}=0}}\\
        &\leq (m_2-2)\!\int_0^{k+1}\! \E{z}{{\rm e}^{\int_0^s\widehat{\xi}(B_r)\, \mathrm{d}r} \widehat{\xi}(B_s) \Bigl(\E{B_s}{{\rm e}^{\int_0^{k+1-s}\widehat{\xi}(B_r)\, \mathrm{d}r} \1{\exists l \in [0\vee k-s,k+1-s) \,:\,B_l=0}}\Bigr)^2 }\, \mathrm{d}s \\
        &\phantom{=} + (m_2-2)\!\int_k^{k+1}\! \E{z}{{\rm e}^{\int_0^s\widehat{\xi}(B_r)\, \mathrm{d}r} \widehat{\xi}(B_s) \1{\exists l \in [k,s)\,:\,B_l=0} \, \Bigl(\E{B_s}{{\rm e}^{\int_0^{k+1-s}\widehat{\xi}(B_r)\, \mathrm{d}r} }\Bigr)^2 }\, \mathrm{d}s;
    \end{split}
    \end{equation}
    here, the first summand corresponds to the case where $\nu_1$ and $\nu_2$ hit the origin \emph{after} the branching, while the second summand corresponds to the complementary case in which $\nu_1$ and $\nu_2$ hit the origin \emph{before} the branching.

    For $s\leq\sqrt{k}$ we can upper bound the inner expectation of the first summand on the right-hand side of \eqref{eq:num_tau_2nd_b0} using Lemma~\ref{lem:num_tau_aux} to obtain  
    \begin{align*}
        \E{B_s}{{\rm e}^{\int_0^{k+1-s}\widehat{\xi}(B_r)\, \mathrm{d}r} \1{\exists l \in [0\vee k-s,k+1-s) \,:\,B_l=0}} \leq (\abs{B_s}\vee1) \frac{C{\rm e}^{\mu_+(k+1-s)}}{(k-s)^{3/2}} {\rm e}^{-\frac{B_s^2}{2(k-s+1)}},
    \end{align*}
    while for $s>\sqrt{k}$, we simply upper bound the indicators by one and $\widehat{\xi}$ by $\mu_+$ to get 
    \begin{align*}
        &\E{z}{{\rm e}^{\int_0^s\widehat{\xi}(B_r)\, \mathrm{d}r} \widehat{\xi}(B_s) \, \Bigl(\E{x}{{\rm e}^{\int_0^{k+1-s}\widehat{\xi}(B_r)\, \mathrm{d}r} \1{\exists l \in [0\vee k-s,k+1-s) \,:\,B_l=0}}_{x=B_s}\Bigr)^2 } \leq C {\rm e}^{\mu_+(2k-s)}.
    \end{align*}
    As a consequence of the previous two estimates, the first integral on the right-hand side of \eqref{eq:num_tau_2nd_b0} can be upper bounded by 
    \begin{align} \label{eq:intermedUBd}
    \begin{split}
        &C{\rm e}^{2\mu_+k}\Bigl( \int_0^{\sqrt{k}} {\rm e}^{-\mu_+ s} \,\E{z}{ \frac{B_s^2\vee1}{(k-s)^3} {\rm e}^{-\frac{B_s^2}{k-s+1}} }\mathrm{d}s + {\rm e}^{-\mu_+ \sqrt{k}}\Bigr)\\
        &\leq C{\rm e}^{2\mu_+k}\Bigl( \int_0^{\sqrt{k}} k^{-3} {\rm e}^{-\mu_+ s}\,\E{z}{(B_s^2+1) {\rm e}^{-\frac{B_s^2}{k-s+1}} }\mathrm{d}s + {\rm e}^{-\mu_+ \sqrt{k}}\Bigr).
        \end{split}
    \end{align}
    In order to further bound this last expression from above, using elementary Gaussian computations we obtain that for any $\lambda>0$,
    \begin{align*}
        \E{z}{{\rm e}^{-\lambda B_s^2}} = \frac{{\rm e}^{-\lambda \frac{z^2}{1+2\lambda s}}}{(1+2\lambda s)^{1/2}} \quad \text{as well as} \quad \E{z}{B_s^2\,{\rm e}^{-\lambda B_s^2}} = \frac{(s + z^2 + 2\lambda s^2) \,{{\rm e}^{-\lambda \frac{z^2}{1+2\lambda s}}}}{(1+2\lambda s)^{5/2}}.
    \end{align*}
    In particular, choosing $\lambda=(k-s+1)^{-1}$, we can continue bounding \eqref{eq:intermedUBd} from above by
    \begin{align*}
        &C{\rm e}^{2\mu_+k}\Bigl(k^{-3} \int_0^{\sqrt{k}} {\rm e}^{-\mu_+ s} (s+z^2+c) \,{\rm e}^{-\frac{z^2}{k+1+s}} \,\mathrm{d}s + {\rm e}^{-\mu_+\sqrt{k}}\Bigr)\\
        &\leq C(z^2\vee1)\,{\rm e}^{2\mu_+k} \max \big \{k^{-3} {\rm e}^{-\frac{z^2}{k+\sqrt{k}}}, {\rm e}^{-\mu_+\sqrt{k}} \big \}.
    \end{align*}
    The second summand on the right-hand side of \eqref{eq:num_tau_2nd_b0} is bounded from above by $C{\rm e}^{\mu_+k}$ (since the indicators can be upper bounded by one and $\widehat{\xi}$ by $\mu_+$). Since this is smaller than our estimate for the first summand, we obtain
    \begin{align*}
        &\Ee{z}{\widehat{\xi}}{\sum_{\substack{\nu_1, \nu_2\in N(k+1) \\ \nu_1 \neq \nu_2}} \1{\exists\,l_1,l_2\in[k,k+1)\,:\,X^{\nu_1}_{l_1}=X^{\nu_2}_{l_2}=0}}\\
        &\quad\leq C(z^2\vee1)\,{\rm e}^{2\mu_+k} \max \big \{k^{-3} {\rm e}^{-\frac{z^2}{k+\sqrt{k}}}, {\rm e}^{-\mu_+\sqrt{k}} \big \}.
    \end{align*}
    Note that this is larger than our estimate of the first moment in \eqref{eq:num_tau_2nd_firstmom}. Going back to \eqref{eq:num_tau_2nd_moments} then finishes the proof.
\end{proof}
We can now use the previous two results to establish that the assumptions of  Lemma~\ref{lem:paley} are indeed satisfied in our context, at least if $z$ is not too large.
We recall the definition of $r^*$ from \eqref{def:r*} as well as the definition of $Z'_t$ in \eqref{eq:Z'}.
\begin{lemma}\label{lem:paley2}
	There exists a constant $\Cl[c]{cst:paley} > 0$ such that for all $z$ with $0\leq z^2 \leq (\sqrt{t}\vee r^*t)^{3/2}$,
	\begin{align*}
        \mathbb{P}_{z}^{\widehat{\xi}}\Bigl(Z^\prime_t(y)\geq\frac{1}{2} \mathbb{E}_{z}^{\widehat{\xi}}[Z^\prime_t(y)]\Bigr) \geq \Cr{cst:paley}
	\end{align*}
    for all $t$ large enough.
\end{lemma}
\begin{proof}
    We use Lemma~\ref{lem:paley} and verify its condition \eqref{cond:paley}.
    By Lemma~\ref{lem:num_tau_lb}, for $k\geq2$, we have the lower bound
    \begin{equation*}
        \left( \Ee{z}{\widehat{\xi}}{\#\{i:\tau^\prime_i \in [k,k+1)\}} \right)^2 \geq C(z\vee1)^2 {\rm e}^{-\frac{z^2}{k}} \frac{{\rm e}^{2\mu_+k}}{k^3},
    \end{equation*}
    while, by Lemma~\ref{lem:num_tau_2nd}, we may upper bound
    \begin{equation*}
        \Ee{z}{\widehat{\xi}}{\bigl(\#\{i:\tau_i^\prime \in [k,k+1)\}\bigr)^2} \leq \frac{C(z\vee 1)^2 {\rm e}^{2\mu_+ k}}{k^3} \max \big\{{\rm e}^{-\frac{z^2}{k+\sqrt{k}}}, k^3 {\rm e}^{-\mu_+ \sqrt{k}}\big\}.
    \end{equation*}
    Consequently,
    \begin{align*}
        \frac{\Ee{z}{\widehat{\xi}}{\bigl(\#\{i:\tau_i^\prime \in [k,k+1)\}\bigr)^2}}{\left( \Ee{z}{\widehat{\xi}}{\#\{i:\tau^\prime_i \in [k,k+1)\}} \right)^2} &\leq C \max
        \big \{{\rm e}^{-\frac{z^2}{k+\sqrt{k}}+\frac{z^2}{k}}, k^3 {\rm e}^{-\mu_+ \sqrt{k}+\frac{z^2}{k}}\big \}\\
        &\leq C \max \big \{{\rm e}^{\frac{z^2}{k^{3/2}}}, k^3 {\rm e}^{\sqrt{k}(\frac{z^2}{k^{3/2}}-\mu_+)} \big\}.
    \end{align*}
	Due to the definition of $\Lambda$ in \eqref{eq:IDef} as well as our assumptions on $z$ we furthermore obtain, for $k \in \Lambda$, 
	\[
    {\frac{z^2}{k^{3/2}}} \leq {\frac{z^2}{\bigl(\sqrt{t}\vee (r^*t-\sqrt{t})\bigr)^{3/2}}} \leq {3/2} < \mu_+ \quad \text{ for all $t$ large enough}.
    \]
    As a consequence of the above, \eqref{cond:paley} is satisfied and Lemma \ref{lem:paley2} follows via Lemma~\ref{lem:paley}.
\end{proof}
\subsubsection{Upper bounding $\Eeb{z}{\widehat{\xi}}{Z_t(y)}$}
In this subsection we provide an upper bound for $\Eeb{z}{\widehat{\xi}}{Z_t(y-\widehat{m}(t))}$ which will then imply an upper bound for $\Pbe{z}{\widehat{\xi}}{M_t \leq y-\widehat{m}(t)}$ via Lemma~\ref{lem:prob_ub}. For this purpose, recall $r^*$ from \eqref{def:r*}, the definition of $\widehat{m}(t)$ in \eqref{def:m(t)},  as well as the definition of $Z_t(y)$ from \eqref{eq:Z}.
\begin{proposition}\label{prop:E_X}
    \begin{enumerate}[label=(\alph*)]
        \item There exists a constant $C>0$ such that for every $\eta_0\in (0,1/2)$, if $y=y(t)=o(\log t)$ (without restriction on the sign) and $z=z(t)\geq0$, then
    \begin{equation}\label{eq:E_X_1}
        \Eeb{z}{\widehat{\xi}}{Z_t(y-\widehat{m}(t))} \leq C (z\vee1) {\rm e}^{\sqrt{2}y} \exp \Big\{-\frac{z^2}{2(r^*t+\sqrt{t}\log t+1)}\Big\} + o(1)
    \end{equation}
    for all sufficiently large $t$, where $o(1)$ denotes a function converging to $0$ as $t\to\infty$, independent of the choice of $z$.

    \item If $\eta_0\in[1/2,\infty]$, there exists a constant $C>0$ such that if $y=y(t)=o(\log t)$ and $z=z(t)\geq0$, then
    \begin{equation}\label{eq:E_X_2.1}
        \Eeb{z}{\widehat{\xi}}{Z_t(y-\widehat{m}(t))} \leq C(z\vee1) {\rm e}^{\sqrt{2}y}\exp\Bigl\{-\frac{z^2}{4\sqrt{t}\log t}\Bigr\} + o(1),
    \end{equation}
    for all $t$ large enough, where $o(1)$ is independent of the choice of $z$.
    \end{enumerate}
\end{proposition}
In order to prove this proposition, we start by upper bounding $\Eeb{z}{\widehat{\xi}}{Z_t(y-\widehat{m}(t))}$ via a deterministic sum. 
\begin{lemma} \label{lem:E_X_int}
    There exists a constant $C>0$ such that,
	with the same notation as in Proposition~\ref{prop:E_X}, we have
	\begin{equation}\label{eq:E_X_int}
		\Eeb{z}{\widehat{\xi}}{Z_t(y-\widehat{m}(t))} \leq C \sum_{k=0}^{\floor{t}} \frac{(z\vee1)}{(k\vee 1)^{3/2}} {\rm e}^{-\frac{z^2}{2(k+1)}} {\rm e}^{\mu_+ k} \, \Pbe{0}{1}{M_{t-k} \leq y - \widehat{m}(t)}.
	\end{equation}
\end{lemma}
\begin{proof}
	By the definition of $Z_t$ in \eqref{eq:Z},
	\begin{equation*}
		\Eeb{z}{\widehat{\xi}}{Z_t(y-\widehat{m}(t))} = \sum_{k=0}^{\floor{t}} \Ee{z}{\widehat{\xi}}{\#\{i\geq1\,:\,\tau_i \in [k,k+1)\}} \, \Pbe{0}{1}{M_{t-k} \leq y -\widehat{m}(t)},
    \end{equation*}
    and the claim follows from Lemma~\ref{lem:num_tau_ub}.
\end{proof}

The next two lemmas state that the sum in \eqref{eq:E_X_int} is concentrated around $k\approx r^*t$ (at least if $\eta_0\in(0,1/2)$). In their proofs we will repeatedly use the elementary bound
\begin{equation}\label{eq:sup_gauss}
    (z\vee1)\,{\rm e}^{-\frac{z^2}{2u}} \leq \sqrt{u}, \qquad z\geq0,\ u\geq1.
\end{equation}
\begin{lemma}\label{lem:E_X_int2.0}
    With the same notation as in Proposition~\ref{prop:E_X}, we have
    \begin{equation}\label{eq:E_X_int2.up}
			\lim_{t\to\infty}\sum_{r^*t+\sqrt{t}\log t<k\leq \floor{t}} \frac{(z\vee1)}{(k\vee1)^{3/2}}\,{\rm e}^{-\frac{z^2}{2(k+1)}} {\rm e}^{\mu_+ k}\,\Pbe{0}{1}{M_{t-k}\leq y-\widehat{m}(t)} = 0
		\end{equation}
    uniformly in $z\geq0$ and every $y=o(\log t)$.
    If $\eta_0\in(0,1/2)$, then $\lim_{t\to\infty}r^*t-\sqrt{t}\log t=\infty$ and we also have
    \begin{equation}\label{eq:E_X_int2.low}
			\lim_{t\to\infty}\sum_{0\leq k<r^*t-\sqrt{t}\log t} \frac{(z\vee1)}{(k\vee1)^{3/2}}\,{\rm e}^{-\frac{z^2}{2(k+1)}} {\rm e}^{\mu_+ k}\,\Pbe{0}{1}{M_{t-k}\leq y-\widehat{m}(t)} = 0
		\end{equation}
    uniformly in $z\geq0$ and every $y=o(\log t)$.
\end{lemma}
Similar to $\Lambda$ in \eqref{eq:IDef}, let
\begin{equation*}
    \Lambda_+ \coloneqq \big[(r^*t-\sqrt{t}\log t)\vee0,\, r^*t+\sqrt{t}\log t\big]\cap\N_0.
\end{equation*}
\begin{lemma}\label{lem:E_X_int2}
	With the same notation as in Proposition~\ref{prop:E_X}, there exists a constant $C>0$ such that
    \begin{align*}
        &\sum_{k\in\Lambda_+} \frac{(z\vee1)}{(k\vee1)^{3/2}}\,{\rm e}^{-\frac{z^2}{2(k+1)}} {\rm e}^{\mu_+ k}\,\Pbe{0}{1}{M_{t-k}\leq y-\widehat{m}(t)}\\
        &\quad\leq C(z\vee1) \sum_{k\in\Lambda_+} \Psi(k, t)\,{\rm e}^{-\frac{z^2}{2(k+1)}} {\rm e}^{-\frac{(k-r^*t)^2}{t}},
    \end{align*}
	where 
    \begin{equation} \label{eq:Psi}
        \Psi(k, t)\coloneqq \frac{v_{+}t - \sqrt{2}(t-k) + \ell_y(t)}{(k\vee1)^{3/2}}\, \exp\Bigl\{- \frac{v_{+}t\, \ell_y(t)}{t-k}\Bigr\}
    \end{equation} 
    with 
    \begin{equation} \label{def:v+}
        v_{+} \coloneqq \mu_+/\sqrt{2(\mu_+-1)}
    \end{equation}
    the prefactor of the linear order term in $\widehat{m}(t)$ (cf.\ \eqref{def:m(t)}), and the logarithmic term
    \begin{equation}\label{eq:ell_y}
    \ell_y(t) \coloneqq \widehat{m}(t)-v_{+}t + \frac{3}{2\sqrt{2}}\log(t-k)-y.
    \end{equation}
    Consequently, when plugged into \eqref{eq:E_X_int}, we obtain
    \begin{equation}\label{eq:E_X_int2}
			\Eeb{z}{\widehat{\xi}}{Z_t(y-\widehat{m}(t))} \leq  C (z\vee1) \sum_{k\in\Lambda_+} \Psi(k, t)\,{\rm e}^{-\frac{z^2}{2(k+1)}} {\rm e}^{-\frac{(k-r^*t)^2}{t}} + o(1),
		\end{equation}
    where the $o(1)$  is the sum of the asymptotically vanishing terms in \eqref{eq:E_X_int2.up} and \eqref{eq:E_X_int2.low}, and is uniform in $z\geq0$.
\end{lemma}
\begin{remark}
    Intuitively, the summation variable $k$ in \eqref{eq:E_X_int} represents the time at which particles leave the region of higher branching rate. These results therefore further support the intuition that $r^*t$ is the optimal exit time for minimal particles from the right half-line, as illustrated in Figure~\ref{fig:optimal_path}.
\end{remark}

\begin{proof}[Proof of Lemma~\ref{lem:E_X_int2.0}]
    Let $k\leq\floor{t}$. By the many-to-one lemma (Lemma \ref{lem:many-to-one}) we have $\Pbe{0}{1}{M_{t-k}\leq y - \widehat{m}(t)} \leq {\rm e}^{t-k} \P_0(B_{t-k}\leq y - \widehat{m}(t))$, and we recall the standard Gaussian estimate (cf.\ (1.2.2) in \cite{AdTa-07}) 
    \begin{align}\label{eq:std_gaussian_ineq}
			\Pb{0}{B_t \leq x} \leq \frac{\sqrt{t}}{\sqrt{2\pi} \abs{x}} {\rm e}^{ -\frac{x^2}{2t} }
		\end{align}
    for $x< 0$. Since also ${\rm e}^{-\frac{z^2}{2(k+1)}}\leq{\rm e}^{-\frac{z^2}{2(t+1)}}$, we obtain
    \begin{equation}\label{eq:E_X_ineq1}
    \begin{split}
        &\frac{(z\vee1)}{(k\vee1)^{3/2}}{\rm e}^{-\frac{z^2}{2(k+1)}} {\rm e}^{\mu_+ k} \, \Pbe{0}{1}{M_{t-k}\leq y - \widehat{m}(t)}\\ 
        &\leq (z\vee1)\,{\rm e}^{-\frac{z^2}{2(t+1)}} \frac{1}{(k\vee1)^{3/2}} \frac{\sqrt{t-k}}{\sqrt{2\pi}(\widehat{m}(t)-y)} {\rm e}^{\mu_+ t - (\mu_+ -1)(t-k) - \frac{(\widehat{m}(t)-y)^2}{2(t-k)}}.
    \end{split}
    \end{equation}
    In order to upper bound \eqref{eq:E_X_ineq1}, we first decompose the exponent on the right-hand side into polynomial and sub-polynomial terms: recalling $\widehat{m}(t)$ from \eqref{def:m(t)} we can write $\widehat{m}(t) - y = v_{+}t -\mathfrak{l}_y(t)$, where 
    \begin{equation*}
    \mathfrak{l}_y(t) \coloneqq \frac{3}{2\sqrt{2}}\log t + \frac{1}{2\sqrt{2}}\log(t^{-\eta}\vee t^{-1/2}) + y;
    \end{equation*}
    note that $0<\mathfrak{l}_y(t)\leq C\log t$ for $t$ large enough, since $y=o(\log t)$.
    We first show that the exponential term on the right-hand side of \eqref{eq:E_X_ineq1} converges to $0$ as $t\to\infty$ uniformly for all $k$ in the range of summation. We start by rewriting
    \begin{equation} \label{eq:expMax}
    \begin{split}
        &\mu_+ t - (\mu_+ -1)(t-k) - \frac{(\widehat{m}(t)-y)^2}{2(t-k)}\\ 
        &= \mu_+ t - (\mu_+ -1)(t-k) - \frac{\bigl(v_{+}t\bigr)^2}{2(t-k)} + \frac{v_{+}t\, \mathfrak{l}_y(t)}{t-k} - \frac{\bigl(\mathfrak{l}_y(t)\bigr)^2}{2(t-k)}.
    \end{split}
    \end{equation}
    Since $v_+=\sqrt{2(\mu_+-1)}\,(1-r^*)$, expanding 
    $$(v_+t)^2=2(\mu_+-1)\,(k-r^*t+t-k)^2$$
    together with the definition of $r^*$ in \eqref{def:r*}, gives
    \begin{equation}\label{eq:E_X_exp_term}
        \mu_+ t - (\mu_+ -1)(t-k) - \frac{\bigl(v_{+}t\bigr)^2}{2(t-k)} = -(\mu_+ -1)\frac{(k-r^*t)^2}{t-k}.
    \end{equation}
   Using \eqref{eq:E_X_exp_term}, the term in \eqref{eq:expMax} is bounded from above by
    \begin{align*}
        \frac{v_{+}t\, \mathfrak{l}_y(t) - (\mu_+ -1)(k-r^*t)^2}{t-k}.
    \end{align*}
    Now note that for $k$ such that $\abs{k-r^*t}>\sqrt{t}\log t$ the numerator satisfies
    \begin{align*}
        v_{+}t\, \mathfrak{l}_y(t) - (\mu_+ -1)(k-r^*t)^2 &\leq v_{+}t\, \mathfrak{l}_y(t) - (\mu_+ -1)t\,(\log t)^2 < 0
    \end{align*}
    for $t$ large enough, uniformly in $k$. Since $0<t-k\leq t$, this implies that uniformly for such $k$, \eqref{eq:expMax} is bounded from above by $v_+\mathfrak{l}_y(t)-(\mu_+-1)(\log t)^2\leq C\log t-(\log t)^2$.
    Furthermore upper bounding $\sqrt{t-k}/(\widehat{m}(t)-y) \leq C/\sqrt{t}$ and $(z\vee1){\rm e}^{-\frac{z^2}{2(t+1)}}\leq\sqrt{t+1}$ by \eqref{eq:sup_gauss} in combination with \eqref{eq:E_X_ineq1}, we thus have
    \begin{equation*}
    \begin{split}
        &\sum_{r^*t+\sqrt{t}\log t<k\leq \floor{t}} \frac{(z\vee1)}{(k\vee1)^{3/2}}{\rm e}^{-\frac{z^2}{2(k+1)}} {\rm e}^{\mu_+ k} \, \Pbe{0}{1}{M_{t-k}\leq y - \widehat{m}(t)} \\
        &\leq C \sum_{k\geq1} \frac{1}{k^{3/2}} \, {\rm e}^{C\log t-(\log t)^2} \longrightarrow 0
    \end{split}
    \end{equation*}
    as $t\to\infty$ uniformly in $z \ge 0$. If $\eta_0\in(0,1/2)$, then $\lim_{t\to\infty}r^*t-\sqrt{t}\log t=\infty$ and the very same computation applies to the summation range $0\leq k<r^*t-\sqrt{t}\log t$, which yields \eqref{eq:E_X_int2.low}.
\end{proof}

\begin{proof}[Proof of Lemma~\ref{lem:E_X_int2}]
    We only need to prove the upper bound
    \begin{equation}\label{eq:E_X_int2_goal}
    \begin{split}
        {\rm e}^{\mu_+ k}\,&\Pbe{0}{1}{M_{t-k}\leq y-\widehat{m}(t)}\\
        &\leq C\bigl(v_{+}t - \sqrt{2}(t-k) + \ell_y(t)\bigr)\, \exp\Bigl\{- \frac{v_{+}t\, \ell_y(t)}{t-k}\Bigr\} \exp\Bigl\{-\frac{(k-r^*t)^2}{t}\Bigr\}
    \end{split}
    \end{equation}
    for all $k\in\Lambda_+$.
    The main difference between this proof and that of Lemma~\ref{lem:E_X_int2.0} is that, by restricting to $k\in\Lambda_+$, we can use \eqref{eq:BramsonAsymp} to upper bound
    \begin{equation}\label{eq:E_X_ineq2}
        {\rm e}^{\mu_+ k}\,\Pbe{0}{1}{M_{t-k}\leq y - \widehat{m}(t)} \leq C x\, {\rm e}^{\mu_+k-\sqrt{2}x - \frac{x^2}{2(t-k)}}
    \end{equation}
    with $x = \widehat{m}(t) - \sqrt{2}(t-k) + \frac{3}{2\sqrt{2}}\log(t-k) - y$. Recalling the notation from \eqref{eq:ell_y} we may write
    \begin{align*}
        x = v_{+}t - \sqrt{2}(t-k) + \ell_y(t).
    \end{align*}
    Here we used that for $k\in\Lambda_+$, $x\geq\ell_y(t)\to\infty$ (cf.\ \eqref{eq:ell_y_eq0} below) and $x=o(t/\log t)$ uniformly in $k$, and that \eqref{eq:BramsonAsymp} holds uniformly over this range of $x$ (cf.\ \cite[Corollary~1.2]{LC25}).
    To show \eqref{eq:E_X_int2_goal}, we need to rewrite the expression in the exponent of \eqref{eq:E_X_ineq2}. We start with
    \begin{align*}
        \frac{x^2}{2(t-k)} &= \frac{\bigl( v_{+}t - \sqrt{2}(t-k) \bigr)^2 + 2\bigl(v_{+}t-\sqrt{2}(t-k)\bigr)\ell_y(t) + \bigl(\ell_y(t)\bigr)^2}{2(t-k)} \\
        &= \frac{\bigl(v_{+}t\bigr)^2}{2(t-k)} -\sqrt{2}v_{+}t + (t-k) +\frac{v_{+}t\, \ell_y(t)}{t-k}- \sqrt{2}\ell_y(t) +  o(1),
    \end{align*}
    and therefore also
    \begin{equation}\label{eq:E_X_expTerm_2}
    \begin{split}
        \mu_+ k -\sqrt{2}x - \frac{x^2}{2(t-k)} &= \mu_+ k + (t-k) - \frac{\bigl(v_{+}t\bigr)^2}{2(t-k)} - \frac{v_{+}t\, \ell_y(t)}{t-k} + o(1) \\
        &= -(\mu_+ -1)\frac{(k-r^*t)^2}{t-k} - \frac{v_{+}t\, \ell_y(t)}{t-k} + o(1)
    \end{split}
    \end{equation}
    by \eqref{eq:E_X_exp_term}, where the $o(1)$ is uniform in $k$.
    Plugging this into \eqref{eq:E_X_ineq2} proves \eqref{eq:E_X_int2_goal}.
\end{proof}

Now that we have analyzed the sum that serves as an upper bound for $\Eeb{z}{\widehat{\xi}}{Z_t(y-\widehat{m}(t))}$ in \eqref{eq:E_X_int}, we are ready to prove Proposition~\ref{prop:E_X}. The exponential factor 
\(
    {\rm e}^{-\frac{(k-r^*t)^2}{t}}
\)
in \eqref{eq:E_X_int2} indicates (at least for small values of $z$) that the principal contribution to the sum comes from indices $k$ around $r^*t$, with fluctuations of order $\sqrt{t}$. However, when $\eta_0\in[1/2,\infty]$, we have $r^*t=O(\sqrt{t})$, so the location of the center is at most of the same order as the fluctuations and therefore does not affect the asymptotic behavior of the sum. In both cases we will use that, since $x\mapsto{\rm e}^{-(x-r^*t)^2/t}$ is unimodal, its sum over the integers is bounded by its maximum plus its integral, i.e.
\begin{equation}\label{eq:disc_gauss}
    \sum_{k\in\Z}{\rm e}^{-\frac{(k-r^*t)^2}{t}}\leq 1+\sqrt{\pi t}.
\end{equation}
\begin{proof}[Proof of Proposition~\ref{prop:E_X}]
	It suffices to upper bound the right-hand side of \eqref{eq:E_X_int2} by \eqref{eq:E_X_1} and \eqref{eq:E_X_2.1}, respectively.
 Recalling $\widehat{m}(t)$ from \eqref{def:m(t)}, we first rewrite $\ell_y(t)$ (cf.\ \eqref{eq:ell_y}) that appears in the function $\Psi(k,t)$ (cf.\ \eqref{eq:Psi}), as
    \begin{equation}\label{eq:ell_y_eq0}
        \ell_y(t) = \frac{3}{2\sqrt{2}}\log\Bigl(1-\frac{k}{t}\Bigr) + \frac{(\eta\wedge\frac12)}{2\sqrt{2}}\log(t)-y = \frac{(\eta\wedge\frac12)}{2\sqrt{2}}\log(t)-y + o(1) \geq 0
    \end{equation}
    for all $t$ large enough, since $k=o(t)$ uniformly for all $k\in\Lambda_+$; note also that the first summand is non-positive, so that $\ell_y(t)\leq \frac{(\eta\wedge\frac12)}{2\sqrt{2}}\log(t)-y$. In particular, using $v_+=\sqrt{2}+O(t^{-2\eta})$ (cf.\ \eqref{eq:taylor}) and $\frac{v_{+}t}{t-k}\geq v_+$, we obtain uniformly for $k\in\Lambda_+$,
    \begin{equation}\label{eq:ell_y_eq1}
        \exp\Bigl\{-\frac{v_{+}t\, \ell_y(t)}{t-k}\Bigr\}\leq {\rm e}^{-v_+\ell_y(t)}= t^{-\frac{\eta\wedge\frac12}{2}}\,{\rm e}^{\sqrt{2}y}\,(1+o(1)).
    \end{equation}
    We distinguish between the two regimes for $\eta_0$.

\bigskip

    \noindent {\em Proof of (a)}: $\eta_0\in(0,1/2)$. In this case, we infer that $r^*t-\sqrt{t}\log t=\Theta(t^{1-\eta})$ for $t$ large enough, so that $\Lambda_+=[r^*t-\sqrt{t}\log t,\, r^*t+\sqrt{t}\log t]\cap\N$, and by \eqref{eq:E_X_int2} we only need to upper bound
    \begin{equation*}
        \sum_{k\in\Lambda_+} \Psi(k, t)\,{\rm e}^{-\frac{z^2}{2(k+1)}} {\rm e}^{-\frac{(k-r^*t)^2}{t}}.
    \end{equation*}
    For all $k\in\Lambda_+$, the definition of $\Psi$ (see \eqref{eq:Psi}) together with \eqref{eq:ell_y_eq1} yields
    \begin{equation}\label{eq:E_X_ineq4}
        \Psi(k, t) \leq \frac{v_{+}t - \sqrt{2}t + \sqrt{2}(r^*t+\sqrt{t}\log t) + \ell_y(t)}{(r^*t-\sqrt{t}\log t)^{3/2}}\, t^{-\frac{\eta}{2}} {\rm e}^{\sqrt{2}y}\,(1+o(1)).
    \end{equation}
    It remains to bound the fraction.
    First we have that
    \begin{align*}
        \frac{v_{+}t - \sqrt{2}t + \ell_y(t)}{(r^*t-\sqrt{t}\log t)^{3/2}} = O\biggl( \frac{(\mu_+ -2)^2 t}{\bigl((\mu_+ -2)t\bigr)^{3/2}} \biggr) = O\big( t^{-\frac{\eta+1}{2}} \big),
    \end{align*}
    using $y=o(\log t)$ and $\eta_0\in(0,1/2)$, and second
    \begin{align*}
        \frac{\sqrt{2}(r^*t+\sqrt{t}\log t)}{(r^*t-\sqrt{t}\log t)^{3/2}}  = \Theta\bigl( t^{\frac{\eta-1}{2}} \bigr),
    \end{align*}
    since $r^*t\pm\sqrt{t}\log t=\Theta(t^{1-\eta})$ for $t$ large enough.
    Therefore, going back to \eqref{eq:E_X_ineq4},
    \begin{equation*}
        \sqrt{t}\,\Psi(k, t) \leq C\, {\rm e}^{\sqrt{2}y}
    \end{equation*}
    uniformly for $k\in\Lambda_+$. Since moreover ${\rm e}^{-\frac{z^2}{2(k+1)}}\leq {\rm e}^{-\frac{z^2}{2(r^*t+\sqrt{t}\log t+1)}}$, Equation \eqref{eq:disc_gauss} yields
    \begin{equation}\label{eq:E_X_ineq5}
    \begin{split}
        \sum_{k\in\Lambda_+} \Psi(k, t)\,{\rm e}^{-\frac{z^2}{2(k+1)}} {\rm e}^{-\frac{(k-r^*t)^2}{t}} &\leq C\, {\rm e}^{-\frac{z^2}{2(r^*t+\sqrt{t}\log t+1)}} {\rm e}^{\sqrt{2}y}\,\frac{1+\sqrt{\pi t}}{\sqrt{t}} \\
        &\leq C\, {\rm e}^{-\frac{z^2}{2(r^*t+\sqrt{t}\log t+1)}} {\rm e}^{\sqrt{2}y}
    \end{split}
    \end{equation}
    for $t$ large enough; together with \eqref{eq:E_X_int2} this proves \eqref{eq:E_X_1} from Proposition~\ref{prop:E_X}.

    \bigskip

    \noindent {\em Proof of (b)}: $\eta_0\in[1/2,\infty]$.
    In this case $r^*t\leq\sqrt{t}$ for $t$ large, so that $\Lambda_+=[0,\, r^*t+\sqrt{t}\log t]\cap\N_0$, and by \eqref{eq:E_X_int2} it is again sufficient to upper bound
    \begin{align}\label{eq:E_X_Case2}
        (z\vee1)\sum_{k\in\Lambda_+} \Psi(k, t)\,{\rm e}^{-\frac{z^2}{2(k+1)}} {\rm e}^{-\frac{(k-r^*t)^2}{t}}
    \end{align}
    accordingly.
    Note that
    \begin{equation*}
    \begin{split}
        \Psi(k, t) &\leq \biggl( \frac{O(1) + \ell_y(t)}{(k\vee1)^{3/2}} + \frac{\sqrt{2}}{(k\vee1)^{1/2}} \biggr) {\rm e}^{- \frac{v_{+}t\, \ell_y(t)}{t-k}} \\
        &\leq t^{-1/4} \biggl( \frac{O(1) + \ell_y(t)}{(k\vee1)^{3/2}} + \frac{\sqrt{2}}{(k\vee1)^{1/2}} \biggr) {\rm e}^{\sqrt{2}y}\,(1+o(1)),
    \end{split}
    \end{equation*}
    using \eqref{eq:Psi}, $v_{+}t-\sqrt{2}t = O((\mu_+-2)^2t) = O(1)$ in this case, and \eqref{eq:ell_y_eq1} with $\eta\wedge\frac12=\frac12$.
    With this bound for $\Psi$ we may upper bound \eqref{eq:E_X_Case2} such that
    \begin{equation}\label{eq:E_X_Case2_int1}
    \begin{split}
        &(z\vee1) \sum_{k\in\Lambda_+} \Psi(k, t)\, {\rm e}^{-\frac{z^2}{2(k+1)}} {\rm e}^{-\frac{(k-r^*t)^2}{t}}\\
        &\leq \frac{O(1) + \ell_y(t)}{t^{1/4}} {\rm e}^{\sqrt{2}y} \sum_{k\in\Lambda_+} \frac{(z\vee1)}{(k\vee1)^{3/2}} {\rm e}^{-\frac{z^2}{2(k+1)}}\,(1+o(1))\\
        &\qquad+ \frac{\sqrt{2}}{t^{1/4}} {\rm e}^{\sqrt{2}y} \sum_{k\in\Lambda_+} \frac{(z\vee1)}{(k\vee1)^{1/2}} {\rm e}^{-\frac{z^2}{2(k+1)}} {\rm e}^{-\frac{(k-r^*t)^2}{t}}\,(1+o(1)).
    \end{split}
    \end{equation}
    By \eqref{eq:sup_gauss}, the summands of the first sum are bounded by $\sqrt{k+1}\,(k\vee1)^{-3/2}\leq\sqrt{2}\,(k\vee1)^{-1}$, so that the first sum is $O(\log t)$ uniformly in $z\geq0$. Since $\ell_y(t)=O(\log t)$ and ${\rm e}^{\sqrt{2}y}=t^{o(1)}$, the first summand on the right-hand side of \eqref{eq:E_X_Case2_int1} converges to $0$ as $t\to\infty$, uniformly in $z\geq0$. Regarding the second sum, we first upper bound ${\rm e}^{-\frac{z^2}{2(k+1)}}\leq{\rm e}^{-\frac{z^2}{4\sqrt{t}\log t}}$, since $r^*t+\sqrt{t}\log t+1 \leq 2\sqrt{t}\log t$ for $t$ large enough.
    Then we split it at $k=\sqrt{t}$ and use \eqref{eq:disc_gauss} for $k>\sqrt{t}$ to obtain
    \begin{equation*}
        \sum_{k\in\Lambda_+} \frac{1}{(k\vee1)^{1/2}}\, {\rm e}^{-\frac{(k-r^*t)^2}{t}} \leq \sum_{0\leq k\leq \sqrt{t}}\frac{1}{(k\vee1)^{1/2}} + \frac{1+\sqrt{\pi t}}{t^{1/4}} \leq C t^{1/4}.
    \end{equation*}
	Together with \eqref{eq:E_X_Case2_int1} this yields
    \begin{equation}\label{eq:E_X_ineq6}
    \begin{split}
        &(z\vee1) \sum_{k\in\Lambda_+} \Psi(k, t)\, {\rm e}^{-\frac{z^2}{2(k+1)}} {\rm e}^{-\frac{(k-r^*t)^2}{t}} \\
        &\leq C(z\vee1){\rm e}^{\sqrt{2}y-\frac{z^2}{4\sqrt{t}\log t}} + o(1)
    \end{split}
    \end{equation}
    where the $o(1)$ terms are uniform in $z$.
    Combined with \eqref{eq:E_X_int2} we infer \eqref{eq:E_X_2.1} from Proposition~\ref{prop:E_X}.
\end{proof}

\subsubsection{Lower bounding $\Eeb{z}{\widehat{\xi}}{Z^\prime_t(y)}$}\label{sec:E_Xprime}
In this subsection, we provide a lower bound for $\Eeb{z}{\widehat{\xi}}{Z^\prime_t(y-\widehat{m}(t))}$ which will then imply a lower bound for $\Pbe{z}{\widehat{\xi}}{M_t \leq y-\widehat{m}(t)}$ via Lemma~\ref{lem:prob_lb}. For this purpose, recall $r^*$ from \eqref{def:r*}, the definition of $\widehat{m}(t)$ in \eqref{def:m(t)}, as well as the definition of $Z'_t(y)$ from \eqref{eq:Z'}.
\begin{proposition}\label{prop:E_Xprime}
    \begin{enumerate}[label=(\alph*)]
    \item There exists a constant $C > 0$ such that for every $\eta_0\in(0,1/2)$, if $y=y(t)=o(\log t)$ and $z=z(t)\geq0$, then 
    \begin{equation*}
        \Eeb{z}{\widehat{\xi}}{Z^\prime_t(y-\widehat{m}(t))} \geq C(z\vee1) \,{\rm e}^{\sqrt{2}y} {\rm e}^{-\frac{z^2}{2(r^*t-\sqrt{t})}}
    \end{equation*}
    for all $t$ sufficiently large.
    
    \item If $\eta_0\in[1/2,\infty]$, there exists a constant $C > 0$ such that, if $y=y(t)=o(\log t)$ and $z=z(t)\geq0$, then
    \begin{equation*}
    		\Eeb{z}{\widehat{\xi}}{Z^\prime_t(y-\widehat{m}(t))} \geq C(z\vee1) \,{\rm e}^{\sqrt{2}y} {\rm e}^{-\frac{z^2}{2\sqrt{t}}}
    \end{equation*}
    for all $t$ large enough.
    \end{enumerate}
\end{proposition}
The proof of Proposition~\ref{prop:E_Xprime} broadly mirrors that of Proposition~\ref{prop:E_X}, with the necessary modifications to derive lower bounds instead of upper bounds.
\begin{lemma}\label{lem:E_Xprime_int2}
		With the same notation as in Proposition~\ref{prop:E_Xprime}, there exists a constant $C\in(0,\infty)$ such that
		\begin{equation}\label{eq:E_Xprime_int2}
			\Ee{z}{\widehat{\xi}}{Z^\prime_t(y-\widehat{m}(t))} \geq  C (z\vee1) \sum_{k\in\Lambda} \Psi(k+1, t)\,{\rm e}^{-\frac{z^2}{2k}} {\rm e}^{-(\mu_+-1)\frac{(k+1-r^*t)^2}{t-k-1}},
		\end{equation}
		with $\Psi$ as in Lemma~\ref{lem:E_X_int2}.
\end{lemma}
\begin{proof}
    Recall the definition of $\Lambda$ (cf.\ \eqref{eq:IDef}) and the definition of $Z^\prime_t$ (cf.\ \eqref{eq:Z'}). By Lemma~\ref{lem:num_tau_lb} we can estimate
		\begin{align}\label{eq:E_Xprime_int}
        \begin{split}
			\Ee{z}{\widehat{\xi}}{Z^\prime_t(y-\widehat{m}(t))} &= \sum_{k\in \Lambda} \Ee{z}{\widehat{\xi}}{\#\{\tau^\prime_i \in [k,k+1)\}} \, \Pbe{0}{1}{M_{t-(k+1)} \leq y -\widehat{m}(t)} \\
			&\geq C \sum_{k\in \Lambda} (z\vee1)\, {\rm e}^{-\frac{z^2}{2k}} {k^{-3/2}} {\rm e}^{\mu_+ k} \, \Pbe{0}{1}{M_{t-(k+1)} \leq y - \widehat{m}(t)}.
        \end{split}
		\end{align}
    For $k\in\Lambda$ we have $k+1\in\Lambda_+$ for $t$ large enough, so, exactly as in the proof of Lemma~\ref{lem:E_X_int2}, but now using the lower bound in \eqref{eq:BramsonAsymp} at time $t-(k+1)$, together with \eqref{eq:E_X_expTerm_2}, we obtain
    \begin{align*}
    		&{\rm e}^{\mu_+ (k+1)}\,\Pbe{0}{1}{M_{t-(k+1)}\leq y-\widehat{m}(t)} \\
            &\geq C \bigl(v_{+}t - \sqrt{2}(t-k-1) + \ell_y(t)\bigr) {\rm e}^{- \frac{v_{+}t\, \ell_y(t)}{t-k-1}} {\rm e}^{-(\mu_+ -1)\frac{(k+1-r^*t)^2}{t-k-1}}\\
        &= C \,(k+1)^{3/2}\,\Psi(k+1, t)\, {\rm e}^{-(\mu_+-1)\frac{(k+1-r^*t)^2}{t-k-1}}.
    \end{align*}
    Plugging this back into \eqref{eq:E_Xprime_int} concludes the proof. 
\end{proof}
Note that the lower bound \eqref{eq:E_Xprime_int2} has the same structure as the upper bound \eqref{eq:E_X_int2}. 
\begin{proof}[Proof of Proposition~\ref{prop:E_Xprime}]
    For all $k\in\Lambda$ we can lower bound
    \begin{equation}\label{eq:E_Xprime_Psi0}
        \Psi(k+1, t) = \frac{v_{+}t - \sqrt{2}(t-k-1) + \ell_y(t)}{(k+1)^{3/2}} {\rm e}^{- \frac{v_{+}t\, \ell_y(t)}{t-k-1}} \geq \frac{\sqrt{2}}{(k+1)^{1/2}} {\rm e}^{- \frac{v_{+}t\, \ell_y(t)}{t-k-1}},
    \end{equation}
    since $v_+\geq\sqrt{2}$ and $\ell_y(t)\geq0$ by \eqref{eq:ell_y_eq0}. We again distinguish the two cases for $\eta_0$.

    \bigskip
    
	\noindent {\em Proof of (a)}: $\eta_0\in(0,1/2)$. In this case we infer that $r^*t-\sqrt{t}=\Theta(t^{1-\eta})$ for $t$ large enough and the set $[r^*t-\sqrt{t},r^*t+\sqrt{t}-1]\cap\N$ is contained in $\Lambda$ and has at least $2\sqrt{t}-2$ elements.
    For $k$ in this set we have ${\rm e}^{-\frac{z^2}{2k}}\geq{\rm e}^{-\frac{z^2}{2(r^*t-\sqrt{t})}}$ and $(k+1)^{-1/2}\geq(r^*t+\sqrt{t})^{-1/2}$; moreover, since $\abs{k+1-r^*t}\leq\sqrt{t}$ and $t-k-1\geq t/2$, we have $(\mu_+-1)\frac{(k+1-r^*t)^2}{t-k-1}=O(1)$.
    Finally, by \eqref{eq:ell_y_eq0}, $\ell_y(t)\leq\frac{\eta}{2\sqrt{2}}\log t-y$, which is non-negative for $t$ large. Hence \eqref{eq:E_Xprime_int2} and \eqref{eq:E_Xprime_Psi0} yield
    \begin{equation}\label{eq:E_Xprime_int3.0}
    		\Ee{z}{\widehat{\xi}}{Z^\prime_t(y-\widehat{m}(t))} \geq  C (z\vee1)\,{\rm e}^{-\frac{z^2}{2(r^*t -\sqrt{t})}}\, \frac{2\sqrt{t}-2}{(r^*t+\sqrt{t})^{1/2}}\, \exp\Bigl\{- \frac{v_{+}t\,(\frac{\eta}{2\sqrt{2}}\log t-y)}{t-r^*t-\sqrt{t}}\Bigr\}.
    \end{equation}
    Since $v_+/(1-r^*-t^{-1/2})=\sqrt{2(\mu_+-1)}\,(1+O(t^{-1/2}))$ and $\sqrt{\mu_+-1}=1+O(t^{-\eta})$, we have
    \begin{align}\label{eq:ell_y_eq3}
        \frac{v_{+}t\,(\frac{\eta}{2\sqrt{2}}\log t-y)}{t-r^*t-\sqrt{t}} = \frac{\eta}{2}\log(t) - \sqrt{2}\, y + o(1).
    \end{align}
    Plugging \eqref{eq:ell_y_eq3} back into \eqref{eq:E_Xprime_int3.0} and using $\sqrt{t}/(r^*t)^{1/2}=\Theta\bigl(t^{\eta/2}\bigr)$ finishes the proof in this case.
    
    \bigskip

    \noindent {\em Proof of (b)}: $\eta_0\in[1/2,\infty]$. In this case $r^*t<\sqrt{t}$ for $t$ large, so the set $[\sqrt{t},2\sqrt{t}-1]\cap\N$ is contained in $\Lambda$ and has at least $\sqrt{t}-1$ elements. For $k$ in this set we have ${\rm e}^{-\frac{z^2}{2k}}\geq{\rm e}^{-\frac{z^2}{2\sqrt{t}}}$, $(k+1)^{-1/2}\geq(2\sqrt{t})^{-1/2}$ and, since $\abs{k+1-r^*t}\leq2\sqrt{t}$, $(\mu_+-1)\frac{(k+1-r^*t)^2}{t-k-1}=O(1)$. Moreover, similarly to \eqref{eq:ell_y_eq3},
    \[
        \frac{v_{+}t\, \ell_y(t)}{t-k-1}\leq\frac{v_{+}t\,(\frac{1}{4\sqrt{2}}\log t-y)}{t-2\sqrt{t}} = \frac{1}{4}\log t-\sqrt{2}y+o(1).
    \]
    Hence \eqref{eq:E_Xprime_int2} and \eqref{eq:E_Xprime_Psi0} yield
    \begin{equation*}
    \begin{split}
        \Ee{z}{\widehat{\xi}}{Z^\prime_t(y-\widehat{m}(t))} &\geq  C (z\vee1)\,{\rm e}^{-\frac{z^2}{2\sqrt{t}}}\,\frac{\sqrt{t}-1}{(2\sqrt{t})^{1/2}}\, t^{-1/4}\,{\rm e}^{\sqrt{2}y}\\
        &\geq C(z\vee1) \,{\rm e}^{\sqrt{2}y} {\rm e}^{-\frac{z^2}{2\sqrt{t}}}
    \end{split}
    \end{equation*}
    for $t$ large enough, which finishes the proof.
\end{proof}

\subsection{Proof of Theorems~\ref{the:min_case1} and \ref{the:min_case2}}\label{sec:proof_bbm_est}
Finally, we can combine the results of this section to prove Theorems~\ref{the:min_case1} and \ref{the:min_case2}.

\begin{proof}[Proof of Theorems~\ref{the:min_case1} and \ref{the:min_case2}]
	The upper bounds in Theorem~\ref{the:min_case1} (Eq.\ \eqref{eq:etaSmallUB}), respectively Theorem~\ref{the:min_case2} (Eq.\ \eqref{eq:etaLargeUB}), follow from Lemma~\ref{lem:prob_ub} together with Proposition~\ref{prop:E_X}.
    Since $z\geq0$, the single term
    \[
        \Pbe{z}{1}{M_t \leq y-\widehat{m}(t)} \leq \Pbe{0}{1}{M_t \leq y-\widehat{m}(t)}
    \]
    tends to $0$ due to \eqref{eq:BramsonAsymp}.
    The lower bounds in Theorem~\ref{the:min_case1} (Eq. \eqref{eq:etaSmallLB}), respectively Theorem~\ref{the:min_case2} (Eq. \eqref{eq:etaLargeLB}), follow from Lemma~\ref{lem:prob_lb} together with Lemma~\ref{lem:paley2} and Proposition~\ref{prop:E_Xprime}. More precisely
	\begin{align}\label{eq:main_proof}
		\Pbe{z}{\widehat{\xi}}{M_t+\widehat{m}(t) \leq y} \geq  \mathbb{P}_{z}^{\widehat{\xi}}\bigl(Z^\prime_t(y-\widehat{m}(t))\geq\frac{1}{2} \mathbb{E}_{z}^{\widehat{\xi}}[Z^\prime_t(y-\widehat{m}(t))]\bigr) \bigl( 1 - {\rm e}^{-\frac{1}{2} \mathbb{E}_{z}^{\widehat{\xi}}[Z^\prime_t(y-\widehat{m}(t))]} \bigr)
	\end{align}
	by Lemma~\ref{lem:prob_lb} and $\mathbb{P}_{z}^{\widehat{\xi}}\bigl(Z^\prime_t(y-\widehat{m}(t))\geq\frac{1}{2} \mathbb{E}_{z}^{\widehat{\xi}}[Z^\prime_t(y-\widehat{m}(t))]\bigr)\geq \Cr{cst:paley}$ by Lemma~\ref{lem:paley2} due to our assumptions on $z$.
    Then, by Proposition~\ref{prop:E_Xprime},
	\begin{align*}
		\mathbb{E}_{z}^{\widehat{\xi}}[Z^\prime_t(y-\widehat{m}(t))] \geq C(z\vee1) \,{\rm e}^{\sqrt{2}y-\frac{z^2}{2((r^*t-\sqrt{t})\vee \sqrt{t})}}
	\end{align*}
	and plugging this into \eqref{eq:main_proof}, the proof of Theorems~\ref{the:min_case1} and \ref{the:min_case2} is finished.
\end{proof}

\subsection{Proof of Theorem~\ref{thm:main_bbm}}\label{sec:proof_bbm_main}
\begin{proof}[Proof of Theorem~\ref{thm:main_bbm}]
     First assume that $\eta_0\in(0,1/2)$. Given $\varepsilon\in(0,1)$, choose $y_-<0$ such that
    \[
        \Cr{cst:ub1} {\rm e}^{\sqrt{2}y_-} < \varepsilon/2.
    \]
    Then, by \eqref{eq:etaSmallUB} of Theorem~\ref{the:min_case1},
    \[
        \Pbe{0}{\widehat{\xi}}{M_t + \widehat{m}(t) \leq y_-} \leq \Cr{cst:ub1}{\rm e}^{\sqrt{2}y_-} + o(1) \leq \varepsilon
    \]
    for all $t$ sufficiently large. On the other hand, choosing $y_+>0$ sufficiently large, \eqref{eq:etaSmallLB} of Theorem~\ref{the:min_case1} implies that
    \[
        \Pbe{0}{\widehat{\xi}}{M_t + \widehat{m}(t) \leq y_+} \geq \Cr{cst:lb11} (1-\varepsilon)
    \]
    for all $t$ sufficiently large.
    In Lemma~\ref{lem:ampl2} below we show that this lower bound can be brought arbitrarily close to $1$ by potentially further increasing the value of $y_+$. Setting $y_0=\max\{y_+,\abs{y_-}\}$, we have shown that there exists $t_0$ such that
    \[
        \sup_{t> t_0}\Pbe{0}{\widehat{\xi}}{\abs{M_t + \widehat{m}(t)} \geq y_0} \leq 2\varepsilon.
    \]
    Tightness for $t\in[1,t_0]$ follows by upper bounding the probability in terms of a homogeneous BBM together with the almost sure boundedness of $M_t + \widehat{m}(t)$. More precisely, for $t\in[1,t_0]$,
    \[
        \Pbe{0}{\widehat{\xi}}{M_t + \widehat{m}(t) \geq y} \leq \Pbe{0}{1}{M_t + \widehat{m}(t) \geq y} \leq \Pbe{0}{1}{\sup_{s\in[1,t_0]} \bigl(M_s + \widehat{m}(s)\bigr) \geq y}
    \]
    and
    \[
        \Pbe{0}{\widehat{\xi}}{M_t + \widehat{m}(t) \leq -y} \leq \Pbe{0}{3}{M_t + \widehat{m}(t) \leq -y} \leq \Pbe{0}{3}{\inf_{s\in[1,t_0]} \bigl(M_s + \widehat{m}(s)\bigr) \leq -y},
    \]
    where we used that $1\leq \widehat{\xi}(z) \leq 3$ for all $z\in\R$. Choosing $y$ sufficiently large then yields tightness for $t\in[1,t_0]$.
    Hence the assertion follows for $\eta_0\in(0,1/2)$.
    
    The case $\eta_0\in[1/2,\infty]$ is treated analogously using Theorem~\ref{the:min_case2}.
\end{proof}

\subsubsection{Amplification}
    In this subsection we prove the amplification result  Lemma~\ref{lem:ampl2} used in the proof of Theorem~\ref{thm:main_bbm} above. We start by recalling a result which states that with high probability we have exponentially many particles close to the origin. Let $N(t,[0,1])$ denote the number of particles located in $[0,1]$ at time $t$.
    \begin{lemma}\label{lem:ampl1}
        For every $(p_k)_{k\in\N}$ fulfilling \eqref{eq:mom_offspr} there exists $\Cl[c]{cst:ampl} = \Cr{cst:ampl}((p_k))>1$ such that, for all $t$ large enough,
        \begin{equation*}
            \Pbe{0}{1}{N(t,[0,1]) \leq \Cr{cst:ampl}^t} \leq \Cr{cst:ampl}^{-t}.
        \end{equation*}
    \end{lemma}
    \begin{proof}
        The assertion follows by a slight adaptation of the proof of Lemma~4.7 in \cite{DreSchmi19}.
    \end{proof}
    \begin{lemma}\label{lem:ampl2}
        Under the assumptions of Theorem~\ref{the:min_case1} or Theorem~\ref{the:min_case2},
        for every $\varepsilon>0$ there exists a $y_0=y_0(\varepsilon)\in\R$ such that
        \begin{equation*}
            \liminf_{t\to\infty}\Pbe{0}{\widehat{\xi}}{M_t + \widehat{m}(t) \leq y} \geq 1-\varepsilon
        \end{equation*}
        for every $y\geq y_0$.
    \end{lemma}
    \begin{proof}
        The proof combines Lemma~\ref{lem:ampl1} with Theorems~\ref{the:min_case1} and \ref{the:min_case2}. Lemma~\ref{lem:ampl1} guarantees that, with probability tending to one, at time $s$ there exist exponentially many particles that are situated close to the origin. For each of these particles, Theorems~\ref{the:min_case1} and \ref{the:min_case2} provide a uniformly positive lower bound on the probability that its descendants contain a minimal particle at the later time $t-s$, provided $y$ is chosen sufficiently large to offset the lost time $s$.

        Restricting to the event $\{ N(s, [0,1]) \geq \Cr{cst:ampl}^s \}$ and conditioning on time $s$, we may estimate
        \begin{equation}\label{eq:ampl-eq1}
        \begin{split}
            &\Pbe{0}{\widehat{\xi}_t}{M_t + \widehat{m}(t) \leq y} \\
            &\geq \Pbe{0}{\widehat{\xi}_{t}}{N(s, [0,1]) \geq \Cr{cst:ampl}^s} \\
            &\phantom{\geq0} \cdot\Bigl\{ 1 - \Ee{0}{\widehat{\xi}_{t}}{ \prod_{\nu\in N(s, [0,1])} \Pbe{X^\nu_s}{\widehat{\xi}_{t}}{M_{t-s} + \widehat{m}(t-s) \geq y-v_{+}s + o(1)} \big\vert N(s, [0,1]) \geq \Cr{cst:ampl}^s} \Bigr\} \\
            &\geq \Pbe{0}{1}{N(s, [0,1]) \geq \Cr{cst:ampl}^s} \Bigl\{ 1 -  \Bigl( \sup_{z\in[0,1]} \Pbe{z}{\widehat{\xi}_{t}}{M_{t-s} + \widehat{m}(t-s) \geq y-v_{+}s + o(1)} \Bigr)^{\Cr{cst:ampl}^s} \Bigr\}
        \end{split}
        \end{equation}
        with $v_{+}$ defined as in \eqref{def:v+}.
        The first factor can be made arbitrarily close to $1$ using Lemma~\ref{lem:ampl1} by taking $s$ large enough. On the other hand, Theorem~\ref{the:min_case1} and Theorem~\ref{the:min_case2} provide us with a universal bound, independent of the starting position, such that
        \begin{equation}\label{eq:ampl-eq2}
        \begin{split}
            &\Bigl( \sup_{z\in[0,1]} \Pbe{z}{\widehat{\xi}_{t}}{M_{t-s} + \widehat{m}(t-s) \geq y-v_{+}s + o(1)} \Bigr)^{\Cr{cst:ampl}^s} \\
            &\leq \Bigl( 1- C \bigl( 1 - \exp\bigl\{-C {\rm e}^{\sqrt{2}(y-v_{+}s + o(1))}\bigr\} \bigr)\Bigr)^{\Cr{cst:ampl}^s}.
        \end{split}
        \end{equation}
        Since $\mu_+(t)-2=(t-s)^{-\eta(t-s)}(1+o(1))$ for fixed $s\geq0$, the proofs of Theorems~\ref{the:min_case1} and \ref{the:min_case2} apply verbatim, with the same constants, to BBM$\widehat\xi_t$ on the time horizon $t-s$, the resulting change of the centering $\widehat m(t-s)$ being $o(1)$.
        Therefore if $y$ is chosen such that $y-v_{+}s\geq s$ for all $t$, then the right hand side of \eqref{eq:ampl-eq2} converges to $0$ by taking first $t\to\infty$ and then $s$ large. Plugging this into \eqref{eq:ampl-eq1} finishes the proof.
    \end{proof}

\subsection{Upper bound when starting below the origin}\label{sec:ub_below_zero}
In this subsection we provide an auxiliary estimate in  case the process starts below zero and cannot immediately profit from the increased branching rate. 
We will use this estimate to upper bound the left $\varepsilon$-front, when proving our main results concerning the diameter of the level sets.
\begin{proposition}\label{prop:ub_below_zero}
	There exists a constant $C>0$ such that for $y=y(t)=o(\log t)$ and $z=z(t)>0$ we have
	\begin{equation*}
		\Pbe{-z}{\widehat{\xi}}{M_t + \widehat{m}(t) \leq y} \leq \Pbe{0}{1}{M_t + \widehat{m}(t)\leq y +z} +  C\,{\rm e}^{\sqrt{2}y - z} + o(1).
	\end{equation*}
\end{proposition}
Taylor's formula implies that there exists a constant $\Cl[c]{cst:asy}>0$ such that
\begin{equation}\label{eq:mhat_lb}
    \widehat{m}(t) - \big(\sqrt{2}t-\frac3{2\sqrt{2}}\log t\big) \geq \Cr{cst:asy} t^{1-2\eta}
\end{equation}
for $t$ sufficiently large.
Together with the asymptotics in \eqref{eq:BramsonAsymp},
Proposition~\ref{prop:ub_below_zero} implies the following corollary.
\begin{corollary}\label{cor:ub_below_zero}
	Let $y\in\R$ be fixed and $z=z(t)$ such that $z\xrightarrow{t\to\infty}\infty$. If 
    \begin{align} \label{eq:zGrowth}
		\Cr{cst:asy} t^{1-2\eta} - z \longrightarrow\infty \quad \text{as} \quad t\to\infty,
	\end{align}
    then
	\begin{equation*}
		\lim_{t\to\infty}\Pbe{-z}{\widehat{\xi}}{M_t+\widehat{m}(t) \leq y} = 0.
	\end{equation*}	
\end{corollary}
\begin{proof}[Proof of Proposition~\ref{prop:ub_below_zero}]
    From Lemma~\ref{lem:prob_ub}, we obtain
    \begin{equation*}
        \Pbe{-z}{\widehat{\xi}}{M_t \leq y} \leq \Pbe{0}{1}{M_t\leq y+z} + \sum_{k=0}^{\floor{t}} \Eeb{-z}{\widehat{\xi}}{\#\{i\geq1\,:\,\tau_i \in [k,k+1)\}} \, \Pbe{0}{1}{M_{t-k} \leq y}.
    \end{equation*}
    Similarly to Lemma~\ref{lem:num_tau_ub} we derive an upper bound for $\Eeb{-z}{\widehat{\xi}}{\#\{i\geq1\,:\,\tau_i \in [k,k+1)\}}$, but this time when starting below the origin; more precisely, we will show below that
    \begin{equation}\label{eq:ub_b0-goal}
        \Eeb{-z}{\widehat{\xi}}{\#\{i\geq1\,:\,\tau_i \in [k,k+1)\}} \leq C\,{\rm e}^{-z} (k\vee1)^{-3/2} {\rm e}^{\mu_+ k}.
    \end{equation}
    Note that this is of the same form as the upper bound in Lemma~\ref{lem:num_tau_ub} when starting from $0$ but with an extra ${\rm e}^{-z}$ term that is independent of $k$. Given \eqref{eq:ub_b0-goal} and arguing as in Lemma~\ref{lem:E_X_int}, we obtain
    \begin{align*}
        &\Pbe{-z}{\widehat{\xi}}{M_t \leq y-\widehat{m}(t)} \\
        &\leq \Pbe{0}{1}{M_t+\widehat{m}(t)\leq y+z} + C\,{\rm e}^{-z} \sum_{k=0}^{\floor{t}} \frac{{\rm e}^{\mu_+ k}}{(k\vee1)^{3/2}} \Pbe{0}{1}{M_{t-k}\leq y-\widehat{m}(t)}.
    \end{align*}
    The sum on the right-hand side coincides with the sum in \eqref{eq:E_X_int} for $z=0$. Following the analysis of this sum in Lemma~\ref{lem:E_X_int2.0}, Lemma~\ref{lem:E_X_int2} and lastly \eqref{eq:E_X_ineq5} as well as \eqref{eq:E_X_ineq6} (with $z=0$),
    we infer
    \begin{equation*}
        \Pbe{-z}{\widehat{\xi}}{M_t + \widehat{m}(t) \leq y} \leq \Pbe{0}{1}{M_t + \widehat{m}(t)\leq y +z} +  C\,{\rm e}^{\sqrt{2}y - z} + o(1)
    \end{equation*}
    and the claim follows. 
    
    It remains to prove \eqref{eq:ub_b0-goal}.
    By an analogous calculation as that leading to \eqref{eq:FKeq}, we have
    \begin{equation*}
        \Eeb{-z}{\widehat{\xi}}{\#\{i\geq1\,:\,\tau_i \in [k,k+1)\}} \leq \E{-z}{{\rm e}^{\int_0^{k+1} \widehat{\xi}(B_r)\, \mathrm{d}r} \1{\exists s\in [k,k+1) \, : \, B_s = 0} }.
    \end{equation*}
    Applying \eqref{eq:num_tau_aux-b2} from Lemma~\ref{lem:num_tau_aux} then yields \eqref{eq:ub_b0-goal}.
\end{proof}

\section{Level sets in a simplified medium}\label{sec:simplePot}
As an intermediate step, we now investigate the diameter of the $\varepsilon$-level set in a simplified medium $\widetilde{\xi}\coloneqq\widetilde{\xi}_t$, for some fixed $t>0$. 
In order to define this medium, recall the definition of $\widehat{\xi}$ in \eqref{def:xiHat} and $\mu_+$ in \eqref{def:mu+}.
For $\varepsilon$ from \eqref{def:front_KPP}, let $\widetilde{x} \coloneqq \widetilde{x}(t, \varepsilon)$ denote the unique real number such that
\begin{equation}\label{def:xTilde}
	\Pbe{0}{\widehat{\xi}}{M_t\leq - \widetilde{x}} = \varepsilon,
\end{equation}
and note that $\widetilde{x}>0$ for all $t$ large enough. Define the simplified medium $\widetilde{\xi}$ by
\begin{equation}\label{def:xiTilde}
	\widetilde{\xi}(z) \coloneqq \begin{cases} 1, & z\in (-\infty, \widetilde{x}), \\ \mu_+, & z\in [\widetilde{x}, \infty), \end{cases}
\end{equation}
so that
\begin{equation}\label{eq:xiTilde_Center}
    \Pbe{\widetilde{x}}{\widetilde{\xi}}{M_t \leq 0} = \varepsilon.
\end{equation}
Similarly to \eqref{def:front_KPP}, and in view of \eqref{eq:mckMin}, we define the right and left $\varepsilon$-fronts in the simplified medium $\widetilde{\xi}$ directly through BBM$\widetilde{\xi}$, i.e.\ via the quantities
\begin{equation}\label{def:mTildePlus}
		\widetilde{m}^{\varepsilon, +}(t) \coloneqq \sup\big\{ x\in\R \,:\, \Pbe{x}{\widetilde{\xi}}{M_t \leq 0} \geq \varepsilon \big\}
	\end{equation}
	and
\begin{equation}\label{def:mTildeMinus}
	\widetilde{m}^{\varepsilon,-}(t) \coloneqq \inf\big\{ x\in\R \,:\, \Pbe{x}{\widetilde{\xi}}{M_t \leq 0} \leq \varepsilon \big\}.
\end{equation}
Observe that by \eqref{eq:xiTilde_Center}, we immediately infer that $\widetilde{m}^{\varepsilon,-}(t) \leq \widetilde{x} \leq \widetilde{m}^{\varepsilon, +}(t)$. 
Analogously to \eqref{def:front_diam}, we define the diameter of the $\varepsilon$-level set in this simplified setting as
$$\widetilde{I}^\varepsilon(t) \coloneqq \widetilde{m}^{\varepsilon, +}(t) - \widetilde{m}^{\varepsilon, -}(t).$$
Recalling that $\widetilde x$ and $\widetilde \xi$ depend on $t$, we then have the following result, which can be considered a first precursor of our main result Theorem \ref{thm:main}; recall $\eta_0$ from \eqref{def:eta}.
\begin{proposition}\label{prop:main_simple}
There exist constants $C>0$ and $\varepsilon_0>0$ such that for each $\varepsilon\in(0,\varepsilon_0)$ fixed, we have:
\begin{enumerate}
    \item 
     For every $\eta_0\in(0,\infty]$,
    \begin{equation}\label{eq:ItildeLB}
        \widetilde{I}^\varepsilon(t) \geq C\,t^{1/4} (\log t)^{1/2}
    \end{equation}
     for all $t$ large enough. 

    \item
    If $\eta_0\in(0,1/2)$, then
    \begin{equation}\label{eq:main_simple1}
        \widetilde{I}^\varepsilon(t) = 
        \Theta\bigl(\max\{(t^{1-\eta}\log t)^{1/2}, t^{1-2\eta}\}\bigr) \quad \text{as } t\to\infty.
    \end{equation}
    In particular, choosing
    \begin{equation*}
    \eta =
    \begin{cases}
        1-2\alpha + \frac{\log(\log(t)\vee1)}{\log(t)}, & \alpha\in(1/4,\frac13], \\
        (1-\alpha)/2, & \alpha\in(\frac13,1),        
    \end{cases}
    \end{equation*}
    yields
    \begin{equation*}
        \widetilde{I}^\varepsilon(t) = 
        \Theta\bigl(t^{\alpha}\bigr) \quad \text{as } t\to\infty
    \end{equation*}
    for all $\alpha\in(1/4,1)$.
\end{enumerate}
\end{proposition}

\begin{remark} \label{rem:sat}
The lower bound in \eqref{eq:ItildeLB} reflects a saturation
phenomenon at the scale
\[
    \mu_+ - 2 = \Theta(t^{-1/2}).
\]
Recall that the optimal particle spends approximately
$r^*t = \Theta(t^{1-\eta})$ units of time in the favorable half-line.
If $\eta_0 \geq 1/2$, then $r^*t = O(\sqrt{t})$, so this time is no
larger than its natural fluctuation scale. Accordingly,
Theorem~\ref{thm:main_bbm} shows that the leading asymptotics of the
minimum no longer depend on the precise value of $\eta$.

The same saturation appears when the process is started at a distance
$z>0$ to the right of the interface. In the regime $\eta_0 \geq 1/2$,
the lower bound in Theorem~\ref{the:min_case2} contains the Gaussian
penalty
\[
    z\exp\left\{-\frac{z^2}{2\sqrt{t}}\right\}.
\]
Consequently, $z$ may be chosen of order
$t^{1/4}\sqrt{\log t}$, with a sufficiently small multiplicative
constant, while the probability of reaching the relevant minimum
remains bounded away from zero. Since such a displacement of the
starting point produces a corresponding displacement of the right
$\varepsilon$-front, this yields
\[
    \widetilde{I}^{\varepsilon}(t)
    \geq C\,t^{1/4}\sqrt{\log t}.
\]
This explains both why the lower bound in
\eqref{eq:ItildeLB} is independent of $\eta$ in the fast-convergence
regime and why varying $\eta$ in the present construction cannot
produce a level set diameter of smaller polynomial order than
$t^{1/4}$.
\end{remark}

\subsection{Proof of Proposition~\ref{prop:main_simple}}
From the tightness result Theorem~\ref{thm:main_bbm} we can infer that for any $\varepsilon\in(0,1)$, there exist $y_-=y_-(\varepsilon)$ and $y_+=y_+(\varepsilon)\in\R$ such that
\begin{equation*}
    \Pbe{0}{\widehat{\xi}}{M_t \leq -\widehat{m}(t) + y} \leq \varepsilon
\end{equation*}
for all $y\leq y_-$, and
\begin{equation*}
    \Pbe{0}{\widehat{\xi}}{M_t \leq -\widehat{m}(t) + y} \geq \varepsilon
\end{equation*}
for all $y\geq y_+$,  for all $t$ large enough. Thus, we may bound
\begin{equation*}
    \widehat{m}(t) - y_+ \leq \widetilde{x} \leq \widehat{m}(t) - y_-
\end{equation*}
such that
\begin{equation}\label{eq:xTilde_Asym}
    \widetilde{x} = \widehat{m}(t) + O(1).
\end{equation}
In order to estimate $\widetilde{m}^{\varepsilon,+}(t)$ and $\widetilde{m}^{\varepsilon,-}(t)$, we compare initial positions relative to $\widetilde{x}$ and then shift the medium so that the process is again viewed with respect to $\widehat{\xi}$. More precisely, for any $z\in\R$ we have
\begin{equation}\label{eq:simple_fronts_eq0}
    \Pbe{\widetilde{x}+z}{\widetilde{\xi}}{M_t \leq 0} = \Pbe{z}{\widehat{\xi}}{M_t \leq -\widetilde{x}} = \Pbe{z}{\widehat{\xi}}{M_t \leq -\widehat{m}(t) +O(1)}
\end{equation}
where the final identity follows from \eqref{eq:xTilde_Asym}. We then apply Theorem~\ref{the:min_case1} and Theorem~\ref{the:min_case2} to select admissible values of $z$ for which the probability in \eqref{eq:simple_fronts_eq0} is either large or small, respectively.
\begin{proof}[Proof of Proposition~\ref{prop:main_simple}]
    {\em Proof of part (a)}: Let $\eta_0\in[1/2,\infty]$.
    Recall the definitions of $\widetilde{m}^{\varepsilon,+}(t)$ in \eqref{def:mTildePlus} and $\widetilde{m}^{\varepsilon,-}(t)$ in \eqref{def:mTildeMinus}. In order to lower bound $\widetilde{m}^{\varepsilon,+}(t)$, we apply \eqref{eq:simple_fronts_eq0} to infer from \eqref{eq:etaLargeLB} of Theorem~\ref{the:min_case2} that there exists a constant $c>0$ such that
    \begin{equation*}
        \Pbe{\widetilde{x}+z}{\widetilde{\xi}}{M_t \leq 0} \geq \Cr{cst:lb21} (1+o(1))
    \end{equation*}
    for $z = (c \sqrt{t}\log t)^{1/2}$ and $t$ sufficiently large. Consequently, for $\varepsilon<\Cr{cst:lb21}$, we have
    \begin{equation*}
        \widetilde{m}^{\varepsilon,+}(t) \geq \widetilde{x} + (c \sqrt{t}\log t)^{1/2}.
    \end{equation*}
    Together with the trivial bound
    \( \widetilde{m}^{\varepsilon,-}(t) \leq \widetilde{x}
    \)
    the previous display yields \eqref{eq:ItildeLB} for $\eta_0\in[1/2,\infty]$. For $\eta_0\in(0,1/2)$, Equation~\eqref{eq:ItildeLB} is an immediate consequence of \eqref{eq:main_simple1} from part (b).

    \bigskip
    
	\noindent {\em Proof of part (b)}: Let $\eta_0\in(0,1/2)$. We are now in the regime in which we can make use of Theorem~\ref{the:min_case1}. Together with \eqref{eq:simple_fronts_eq0}, for $\varepsilon<\Cr{cst:lb11}$ we obtain the bounds
    \begin{equation*}
        \widetilde{x} + c\bigl((r^*t)\log t\bigr)^{1/2}  \leq \widetilde{m}^{\varepsilon,+}(t) \leq \widetilde{x} + C\bigl((r^*t)\log t\bigr)^{1/2},
    \end{equation*}
    for some constants $c,C>0$.
    Since $r^*t = \Theta\bigl(t^{1-\eta}\bigr)$ this implies
    \begin{equation}\label{eq:simple_rf-1}
        \widetilde{m}^{\varepsilon,+}(t) - \widetilde{x} = \Theta\bigl((t^{1-\eta}\log t)^{1/2}\bigr).
    \end{equation}
    In order to understand $\widetilde{m}^{\varepsilon,-}(t)$, we use the simple observation
    \begin{equation*}
        \Pbe{z}{\widetilde{\xi}}{M_t\leq0} \geq \Pbe{z}{1}{M_t\leq0}
    \end{equation*}
    which---together with Bramson's classical result \cite{bramson1983convergence}
    for constant branching rate $1$---entails that
    \begin{equation*}
        \widetilde{m}^{\varepsilon, -}(t) \geq \sqrt{2}t-\frac{3}{2\sqrt{2}}\log t-C.
    \end{equation*}
    Furthermore, since $\eta_0<1/2$, Taylor's formula (cf. \eqref{eq:taylor}) implies 
    $$
    \widehat{m}(t) - \Big(\sqrt{2}t-\frac{3}{2\sqrt{2}}\log t \Big) = \Theta\bigl((\mu_+-2)^2t\bigr) = \Theta\bigl(t^{1-2\eta}\bigr)$$
    and with $\widetilde{x} = \widehat{m}(t)+O(1)$ by \eqref{eq:xTilde_Asym}, we get
    \begin{equation*}
        \widetilde{m}^{\varepsilon, -}(t) \geq \widetilde{x} - Ct^{1-2\eta}.
    \end{equation*}
    For the upper bound we appeal to Corollary~\ref{cor:ub_below_zero} which, together with \eqref{eq:simple_fronts_eq0}, yields
    \begin{equation*}
        \widetilde{m}^{\varepsilon,-}(t) \leq \widetilde{x} - Ct^{1-2\eta}.
    \end{equation*}
    Hence we infer
    \begin{equation}\label{eq:simple_lf-1}
        \widetilde{x} - \widetilde{m}^{\varepsilon,-}(t) = \Theta\bigl(t^{1-2\eta}\bigr).
    \end{equation}
    Combining \eqref{eq:simple_lf-1} and \eqref{eq:simple_rf-1} then implies \eqref{eq:main_simple1}, and the proof is finished.
\end{proof}

\section{Transfer from simplified medium to full setting}\label{sec:fullPot}
    After considering various auxiliary media, we finally define the medium $\xi$ used to prove Theorems~\ref{thm:main} and \ref{thm:main_crit}. The guiding idea is to construct a medium that consists of increasingly long stretches where the medium alternates between the values $1$ and $\mu_n$, where $\mu_n\to2$ as $n\to\infty$.
    The lengths and locations of these stretches are chosen so that, at suitable times $t_n$, the extremal behavior of BBM$\xi$ is essentially governed by a single long stretch on which the medium equals $1$ on the left part of the stretch and $\mu_n$ on the right part of the stretch. Later stretches lie too far away to be reached on the time scale $t_n$, whereas earlier stretches are too short to significantly influence $M_{t_n}$.
    
	We start by defining $\xi$ as follows. Choose $t_1 = 1$ and for $t_n$ given define 
    \begin{equation}\label{def:tn}
        t_{n+1} \coloneqq {\rm e}^{({\rm e}^{t_n})}, \quad n \ge 1.
     \end{equation}
     Let
     \begin{equation}\label{def:mu_n}
        \mu_n \coloneqq \mu_+(t_n) = 2 + t_n^{-\eta}
    \end{equation}
    with $\eta=\eta(t_n)$ analogous to \eqref{def:mu+} and \eqref{def:eta}.
    With $y_0=0$, we successively define, given $y_{n-1}$,
    \begin{equation}\label{def:xn}
        x_n \coloneqq x_n(\varepsilon) \coloneqq \widetilde{x}(t_n,\varepsilon/2) + y_{n-1} + 2, \quad n\in\N
    \end{equation}
    and
    \begin{equation}\label{def:y_n}
		y_n \coloneqq x_n + \sqrt{6}t_n,
	\end{equation}
    so that
    \begin{equation}\label{eq:xn_prob}
        \Pbe{0}{\widehat{\xi}}{M_{t_n} \leq -x_n+y_{n-1}+2} = \varepsilon/2 \qquad \forall\,n\in\N,
    \end{equation}
    by choice of $x_n$ in \eqref{def:xn} and the definition of $\widetilde{x}$ in \eqref{def:xTilde}. (This particular choice of $x_n$ is due to its compatibility with our upper bounds in Lemmas~\ref{lem:gen_ub_+} and \ref{lem:gen_ub_-}.)
    Recalling the asymptotics \eqref{eq:xTilde_Asym} we have
	\begin{align}\label{eq:x_n}
		x_n - y_{n-1} = \widehat{m}(t_n) + O(1),
	\end{align}
    for $n$ large enough.
	We define the medium $\xi$ as
	\begin{equation}\label{def:xi}
		\xi(z) \coloneqq \begin{cases} 
		1, & z\in(-\infty, x_1-1)\cup\bigcup_{n\geq 1} (y_n+1, x_{n+1}-1), \\[0.8ex]
		\mu_n, & z\in[ x_n, y_n], \\[0.8ex]
		(x_n - z) + \mu_n (z-(x_n-1)), & z\in[ x_n -1, x_n), \\[0.8ex]
		(z - y_n) + \mu_n (y_n+1 - z), & z\in( y_n, y_n +1]. \\
		\end{cases}
	\end{equation}

	\begin{figure}[htbp]
    	\centering
    	\includegraphics[width=\textwidth]{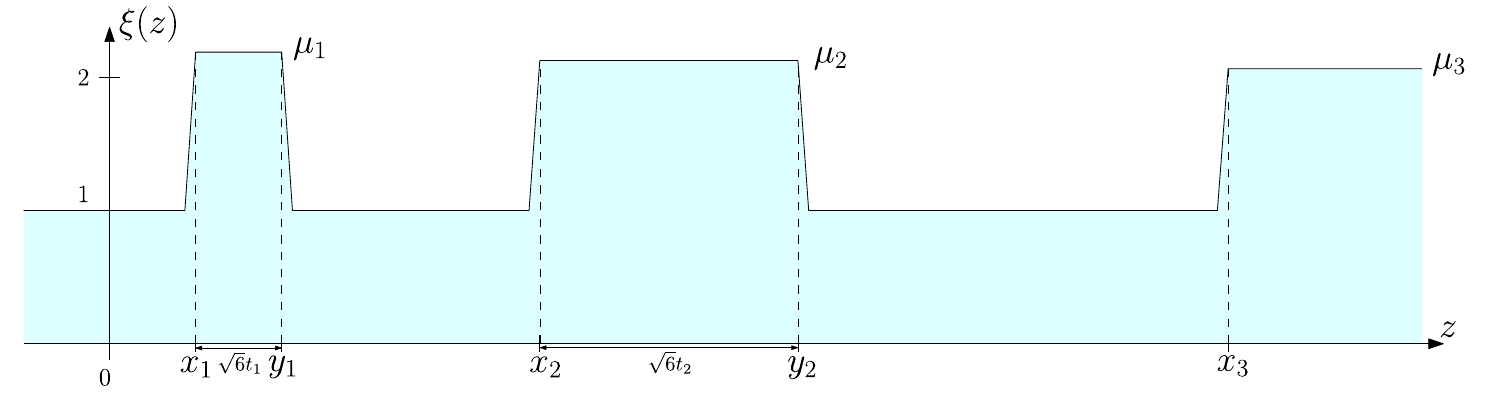}
    	\caption{Sketch of the medium $\xi$.}
    	\label{fig:xi}
	\end{figure}
   
\noindent Figure~\ref{fig:xi} provides a visualization of the medium $\xi$.
Since $\xi$ as defined above is locally H\"{o}lder continuous, we can then use \eqref{eq:McKean} again to get the correspondence of results pertaining to the $\varepsilon$-fronts of solutions to \eqref{eq:KPP}, and those concerning the distributions of the minima of BBM$\xi$  (cf.\ \eqref{def:front_KPP}),  i.e.
	\begin{equation*}
		m^{\varepsilon, +}(t) = \sup\left\{ x\in\R \,:\, u(t,x) \geq \varepsilon \right\} = \sup\{ x\in\R \,:\, \Pbe{x}{\xi}{M_t \leq 0} \geq \varepsilon \}
	\end{equation*}
	and
	\begin{equation*}
		m^{\varepsilon, -}(t) = \inf\{ x\in\R \,:\, u(t,x) \leq \varepsilon \} = \inf\{ x\in\R \,:\, \Pbe{x}{\xi}{M_t \leq 0} \leq \varepsilon \}.
	\end{equation*}
    Recall the definition $I^\varepsilon(t)=m^{\varepsilon, +}(t)-m^{\varepsilon, -}(t)$ from \eqref{def:front_diam}.
    We have the following result:
 \begin{theorem}\label{thm:main_gen}
        In the above setting, there exists a constant $\varepsilon_0>0$ such that for all $\varepsilon\in(0,\varepsilon_0)$ fixed, there exists a constant $C>0$ such that
        \begin{equation}\label{eq:ILB}
            I^\varepsilon(t_n) \geq C\,t_n^{1/4} (\log t_n)^{1/2}
        \end{equation}
        for all $\eta_0\in(0,\infty]$ and all $n$ large enough. On the other hand, if $\eta_0\in(0,1/2)$, then
        \begin{equation}\label{eq:main_1}
            I^\varepsilon(t_n) = 
            \Theta\bigl(\max\{(t_n^{1-\eta}\log t_n)^{1/2}, t_n^{1-2\eta}\}\bigr) \quad \text{as } n\to\infty.
        \end{equation}
    \end{theorem}
    Our main results Theorem~\ref{thm:main} and Theorem~\ref{thm:main_crit} follow almost immediately from Theorem~\ref{thm:main_gen}; see Section~\ref{sec:main_proofs} for the proofs of Theorem~\ref{thm:main_gen} and Theorems~\ref{thm:main} and \ref{thm:main_crit}.
    
    As a preparation for the proof of Theorem~\ref{thm:main_gen} we now show how to transfer the bounds from the setting of the simplified medium to the general setting. Analogously to \eqref{def:r*}, define
    \begin{equation}\label{def:r^*n}
        r^*_n \coloneqq \frac{\mu_n-2}{2(\mu_n-1)} \sim t_n^{-\eta}/2.
    \end{equation}
    
	\subsection{Upper bound for $m^{\varepsilon, +}(t_n)$}
    \begin{lemma}\label{lem:gen_ub_+}
		There exists a constant $C>0$ such that, setting 
        $$z_n = \bigl(C(r^*_n t_n)\log t_n\bigr)^{1/2},$$
        we have for all $\eta_0\in(0,1/2)$,
		\begin{equation*}
			\Pbe{x_n + z}{\xi}{M_{t_n} \leq 0} \leq \varepsilon \qquad \forall z\geq z_n
		\end{equation*}
        for all $n$ large enough.
        In particular, it follows that
		\begin{equation*}
			m^{\varepsilon, +}(t_n) \leq x_n + z_n
		\end{equation*}
        for all $n$ large enough.
	\end{lemma}
    \begin{proof}
    Since it will be easier to work with piecewise constant media, in order to estimate $\xi$ we introduce  $\xi^\prime = \xi'_n$  defined as
	\begin{align*}
		\xi^\prime(z) \coloneqq \begin{cases} \mu_n, & z \geq 0, \\ \mu_1, & z\in[ -x_n+1, -x_n +y_{n-1}+2 ], \\ 1, & \text{otherwise}. \end{cases}
	\end{align*}
 
    \begin{figure}[htbp]
    	\centering
    	\includegraphics[width=\textwidth]{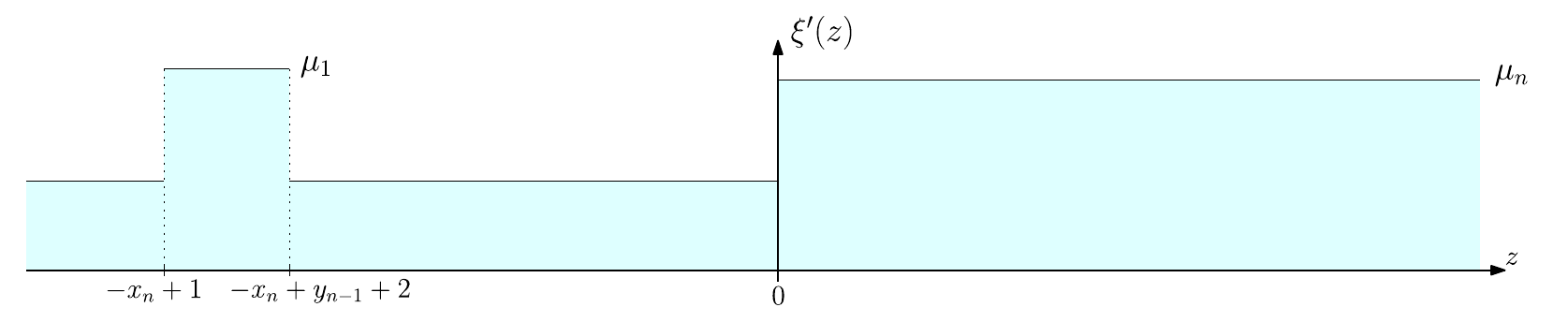}
    	\caption{Sketch of the medium $\xi^\prime$.}
	\end{figure}
    Since $\mu_n$ (cf. \eqref{def:mu_n}) is decreasing in $n$, we have  $\xi(x_n-1 + z) \leq \xi^\prime(z)$ for all $z\in\R$ and we may estimate
	\begin{align}\label{eq:gen_ub_bd0.1}
		\Pbe{x_n+z}{\xi}{M_{t_n} \leq 0} \leq \Pbe{z+1}{\xi^\prime}{M_{t_n} \leq - x_n+1}.
	\end{align}
    Recall the notation of the running minimum in \eqref{eq:runMin}.
	The idea is to replace the event $\left\{ M_{t_n} \leq - x_n+1 \right\}$, which requires the leftmost particle to lie below $-x_n+1$ at time $t_n$, by the superset $\left\{ \underline{M}_{t_n} \leq - x_n + y_{n-1} + 2 \right\}$ which only requires that some particle reaches a position below $-x_n + y_{n-1} +2$ at some time prior to $t_n$; so
    \begin{align*}
        \Pbe{z+1}{\xi^\prime}{M_{t_n} \!\leq \!- x_n+1} \!\leq\! \Pbe{z+1}{\xi^\prime}{\underline{M}_{t_n} \!\leq \!- x_n + y_{n-1} + 2}
    \end{align*}
    since $y_{n-1}\geq0$ (see \eqref{def:y_n}). Moreover,
	\begin{align*}
        \Pbe{z+1}{\xi^\prime}{\underline{M}_{t_n} \!\leq \!- x_n + y_{n-1} + 2} \!=\! \Pbe{z+1}{\widehat{\xi}}{\underline{M}_{t_n} \!\leq \!-x_n + y_{n-1}+2}
	\end{align*}
    since the first hitting of level $-x_n+y_{n-1}+2$ happens by a particle that is born to the right of this location.
    Importantly, observe here that we have now reduced our estimates to probabilities involving the medium $\widehat \xi = \widehat \xi_{t_n}$ introduced in \eqref{def:xiHat}. Now, by Lemma~\ref{lem:BBM_refl}, we may replace $\underline{M}_{t_n}$ with $M_{t_n}$ to obtain
    \begin{equation}\label{eq:gen_ub_bd0}
    		\Pbe{z+1}{\xi^\prime}{M_{t_n} \leq - x_n +1} \leq 2\, \Pbe{z+1}{\widehat{\xi}}{{M}_{t_n} \leq -x_n + y_{n-1}+2}.
    \end{equation}
    Then, replacing $x_n$ by its asymptotics in \eqref{eq:x_n}, together with \eqref{eq:gen_ub_bd0.1} we get that there exists a constant $C>0$ such that
    \begin{equation}\label{eq:gen_ub_b1}
        \Pbe{x_n+z}{\xi}{M_{t_n} \leq 0} \leq 2 \Pbe{z+1}{\widehat{\xi}}{{M}_{t_n} + \widehat{m}(t_n) \leq C}
    \end{equation}
    for all $n$ sufficiently large. Having estimated our quantity of interest in terms of the medium $\widehat{\xi}$ we may now take advantage of Theorem~\ref{the:min_case1} to further bound \eqref{eq:gen_ub_b1}; therefore, if
    \begin{equation}\label{eq:gen_ub_cond}
        (z+1)\,{\rm e}^{-\frac{(z+1)^2}{2r^*_nt_n}} \to 0
    \end{equation}
    as $n\to\infty$, then the claim follows.
    By choosing $z_n = \bigl(C(r^*_n t_n)\log t_n\bigr)^{1/2}$ we obtain
    \begin{align*}
        z_n\,{\rm e}^{-\frac{z_n^2}{2r^*_nt_n}} &= \bigl(C(r^*_nt_n)\log t_n\bigr)^{1/2} t_n^{-C/2}
    \end{align*}
    from which we infer that \eqref{eq:gen_ub_cond} is fulfilled for all $z\geq z_n$ if $C>0$ is chosen sufficiently large.
    \end{proof}

    \subsection{Upper bound for $m^{\varepsilon, -}(t_n)$}
    \begin{lemma}\label{lem:gen_ub_-}
        We have
        \begin{equation}\label{eq:gen_ub_-_mu=2}
            \Pbe{x_n-1}{\xi}{M_{t_n} \leq 0} \leq \varepsilon
        \end{equation}
        for all $n\in\N$.
        Moreover, for any $\eta_0\in(0,1/2)$, setting
        $$z_n = \frac{\Cr{cst:asy}}{2} t^{1-2\eta}_n,$$
        with $\Cr{cst:asy}$ denoting the same constant as in \eqref{eq:mhat_lb},
        we have
        \begin{equation} \label{eq:PMtoZero}
            \lim_{n\to\infty}\Pbe{x_n-z_n}{\xi}{M_{t_n} \leq 0} = 0.
        \end{equation}
        In particular,
        \begin{equation*}
            m^{\varepsilon,-}(t_n) \leq \begin{cases}
                x_n - z_n, & \eta_0\in(0,1/2), \\
                x_n -1, & \eta_0\in[1/2,\infty],
            \end{cases}
        \end{equation*}
        for all $n$ large enough.
    \end{lemma}
    \begin{proof}
        By a calculation analogous to that leading to  \eqref{eq:gen_ub_bd0}, we obtain
        \begin{equation*}
            \Pbe{x_n-z}{\xi}{M_{t_n} \leq 0} \leq 2 \Pbe{-z+1}{\widehat{\xi}}{{M}_{t_n} \leq -x_n + y_{n-1} + 2}
        \end{equation*}
        for any $z>0$. Then \eqref{eq:gen_ub_-_mu=2} follows immediately from \eqref{eq:xn_prob}.

        Let $\eta_0\in(0,1/2)$ from now on. Replacing $x_n$ by its asymptotics in \eqref{eq:x_n} yields 
        \begin{equation*}
            \Pbe{x_n-z_n}{\xi}{M_{t_n} \leq 0} \leq 2 \Pbe{-z_n+1}{\widehat{\xi}}{{M}_{t_n} +\widehat{m}(t_n) \leq C},
        \end{equation*}
        to which Corollary~\ref{cor:ub_below_zero} can now be applied.
        Since $\eta_0<1/2$, we infer
        \[
            z_n = \frac{\Cr{cst:asy}}{2}t_n^{1-2\eta}\xrightarrow{n\to\infty} \infty.
        \]
        Moreover, the growth condition \eqref{eq:zGrowth} for $z_n$ is fulfilled by choice of $z_n$.
        Therefore, the conditions for Corollary~\ref{cor:ub_below_zero} are fulfilled, from which \eqref{eq:PMtoZero} follows.
    \end{proof}
	
	\subsection{Lower bound for $m^{\varepsilon, +}(t_n)$}
	\begin{lemma}\label{lem:RFLB}
        There exists a constant $\varepsilon_0>0$ such that for any $\varepsilon\in(0,\varepsilon_0)$,
		there exists a constant $C>0$ such that setting
        $$z_n= \bigl(C(\sqrt{t_n}\vee(r^*_nt_n))\log t_n\bigr)^{1/2}$$
        we have
		\begin{equation*}
			\liminf_{n\to\infty}\Pbe{x_n+z_n}{\xi}{M_{t_n}\leq 0} > \varepsilon.
		\end{equation*}
		In particular,
		\begin{equation*}
			m^{\varepsilon, +}(t_n) \geq x_n + z_n,
		\end{equation*}
        for all $n$ large enough.
	\end{lemma}
	\begin{proof}
    In contrast to Lemma~\ref{lem:gen_ub_+} we now have to consider that $\xi$ is not constant equal to $\mu_n$ to the right of $x_n$; it is therefore useful to restrict to those particles that stay below the end of the $\mu_n$-stretch. In analogy to \eqref{eq:runMin}, for a stochastic process $(X_t)$ and $t \in (0,\infty)$, we write 
    \begin{equation*}
	    \overline{X}_t := \sup_{s\in[0,t]} X_s
    \end{equation*}
    for the running maximum.
		We start by lower bounding
		\begin{align} \label{eq:withinStretch}
        \begin{split}
			\Pbe{x_n+z_n}{\xi}{M_{t_n}\leq 0} &\geq \Pbe{x_n+z_n}{\widetilde{\xi}^{'}}{M_{t_n} \leq 0, \overline{X}^\nu_{t_n} \leq y_n \text{ for all } \nu\in N(t_n) } \\
			&= \Pbe{x_n+z_n}{\widetilde{\xi}^{'}}{M_{t_n} \leq 0} - \Pbe{x_n+z_n}{\widetilde{\xi}^{'}}{M_{t_n} \leq 0, \,\exists\, \nu \in N(t_n) \text{ s.t. } \overline{X}^\nu_{t_n} > y_n},
        \end{split}
		\end{align}
        where we write $\widetilde{\xi}^{'} = \widetilde{\xi}^{'}_{t_n}$ defined as $\widetilde{\xi}$ in \eqref{def:xiTilde} but centered at $x_n$.
        Note that in the lower bound, we simply disregard all stretches of increased branching rates that preceded the $\mu_n$-stretch.
		  Markov's inequality together with the many-to-one lemma imply
		\begin{align*}
			\Pbe{x_n+z_n}{\widetilde{\xi}^{'}}{\exists\, \nu \in N(t_n) \text{ s.t. } \overline{X}^\nu_{t_n} > y_n}\leq {\rm e}^{\mu_n t_n} \Pb{0}{\overline{B}_{t_n} > y_n-x_n-z_n}.
        \end{align*}
        Applying the reflection principle together with the standard Gaussian estimate \eqref{eq:std_gaussian_ineq}, the right-hand side of the previous display can further be upper bounded by
        \begin{align*}
            & 2 {\rm e}^{\mu_n t_n} \Pb{0}{B_{t_n} > y_n-x_n-z_n}
            \leq 2 \exp\left\{ \mu_n t_n - \frac{(\sqrt{6}t_n-z_n)^2}{2t_n} \right\} \\
            &= 2 \exp\left\{ \mu_n t_n - \frac{6 t_n^2(1-\frac{z_n}{\sqrt{6}t_n})^2}{2t_n} \right\};
		\end{align*}
        this then implies
        \begin{equation}\label{eq:lb_gen_me+_eq1}
            \lim_{n\to\infty}\Pbe{x_n+z_n}{\widetilde{\xi}^{'}}{M_{t_n} \leq 0, \,\exists\, \nu \in N(t_n) \text{ s.t. } \overline{X}^\nu_{t_n} > y_n} = 0,
        \end{equation}
		since $\lim_{n\to\infty}\mu_n = 2$ and $z_n/t_n = o(1)$. Centering the medium back to the origin together with the asymptotics in \eqref{eq:x_n}, we also have
        \begin{equation}\label{eq:lb_gen_me+_eq2}
        \begin{split}
            &\Pbe{x_n+z_n}{\widetilde{\xi}^{'}}{M_{t_n} \leq 0} = \Pbe{z_n}{\widehat{\xi}}{M_{t_n} + x_n \leq 0}\\
            &\geq  \Pbe{z_n}{\widehat{\xi}}{M_{t_n} + \widehat{m}(t_n) \leq -y_{n-1}-C}\geq \min\{\Cr{cst:lb11},\Cr{cst:lb21}\} \, (1+o(1))
        \end{split}
        \end{equation}
        due to Theorems~\ref{the:min_case1} and \ref{the:min_case2} if the constant $C>0$ in
        $$
            z_n=\bigl(C(\sqrt{t_n}\vee(r^*_nt_n))\log t_n\bigr)^{1/2}
        $$
        is chosen sufficiently small  so that
        \[
            z_n\,{\rm e}^{-\sqrt{2}y_{n-1}-\frac{z_n^2}{2(\sqrt{t_n}\vee(r^*_nt_n))}} \to \infty
        \]
        as $n\to\infty$.
        As a consequence, we can plug \eqref{eq:lb_gen_me+_eq1} and \eqref{eq:lb_gen_me+_eq2} back into \eqref{eq:withinStretch} to obtain
  		\begin{align*}
			\Pbe{x_n+z_n}{\xi}{M_{t_n}\leq 0} 
			\geq \min\{\Cr{cst:lb11},\Cr{cst:lb21}\} + o(1) \geq \varepsilon,
		\end{align*}
    for $n$ large enough, which finishes the proof.
	\end{proof}
	
	\subsection{Lower bound for $m^{\varepsilon, -}(t_n)$}
	\begin{lemma}\label{lem:LFLB}
		There exists a constant $C>0$ such that for all $n\in\N$ large enough we have
		\begin{align*}
			m^{\varepsilon, -} (t_n) \geq \sqrt{2}t_n - \frac{3}{2\sqrt{2}} \log t_n - C.
		\end{align*}
	\end{lemma}
	\begin{proof}
		Taking advantage of Bramson's classical result \cite{bramson1983convergence}
  for constant branching rate $1$ in combination with
		\begin{align*}
			\Pbe{z}{\xi}{M_t \leq 0} \geq \Pbe{z}{1}{M_t \leq 0}
		\end{align*}
        already provides us with the desired result.
	\end{proof}
	
	\subsection{Proof of Theorem~\ref{thm:main}, Theorem \ref{thm:main_crit} and Theorem~\ref{thm:main_gen}}\label{sec:main_proofs}
        We first prove Theorem~\ref{thm:main_gen}.
        \begin{proof}[Proof of Theorem~\ref{thm:main_gen}]
            We need to bound $I^\varepsilon(t_n) = m^{\varepsilon, +}(t_n) - m^{\varepsilon, -}(t_n)$.
            Combining Lemma~\ref{lem:gen_ub_-} with Lemma~\ref{lem:RFLB} yields the lower bound
            \begin{equation*}
            \begin{split}
                I^\varepsilon(t_n) &\geq x_n + \bigl(C(\sqrt{t_n}\vee(r^*_nt_n))\log t_n\bigr)^{1/2} - x_n + C\bigl(t_n^{1-2\eta}\vee1\bigr) \\
                &\geq C \max\bigl\{ t_n^{1-2\eta},  \bigl(t_n^{1/4}\vee t_n^{(1-\eta)/2}\bigr)(\log t_n)^{1/2}\bigr\}
            \end{split}
            \end{equation*}
            for all $\eta_0\in(0,\infty]$, from which we can already infer \eqref{eq:ILB}.         
            The upper bound for $\eta_0\in(0,1/2)$ follows from Lemma~\ref{lem:gen_ub_+} and Lemma~\ref{lem:LFLB}:
            \begin{equation*}
            \begin{split}
                I^\varepsilon(t_n) &\leq x_n + \bigl(C(r^*_nt_n)\log(t_n)\bigr)^{1/2} - \bigl( \sqrt{2}t_n - \frac{3}{2\sqrt{2}} \log t_n - C \bigr) \\
                &\leq C\max\bigl\{ t_n^{1-2\eta}, t_n^{(1-\eta)/2}(\log t_n)^{1/2} \bigr\},
            \end{split}
            \end{equation*}
            by \eqref{def:r^*n} and
            \[
                x_n - \big(\sqrt{2}t_n-\frac{3}{2\sqrt{2}}\log t_n\big) = \Theta((\mu_n-2)^2t_n) = \Theta(t_n^{1-2\eta})
            \]
            for $\eta_0\in(0,1/2)$ due to \eqref{eq:x_n} and $y_{n-1}=o(\log t_n)$ by choice of $t_n$ in \eqref{def:tn}. This proves \eqref{eq:main_1} and the proof is finished.
        \end{proof}
        Theorems~\ref{thm:main} and \ref{thm:main_crit} are direct consequences of Theorem~\ref{thm:main_gen}.
        \begin{proof}[Proof of Theorem~\ref{thm:main} and Theorem \ref{thm:main_crit}]
             When choosing
        \begin{equation*}
        \eta =
        \begin{cases}
            1-2\alpha + \frac{\log(\log (t)\vee1)}{\log (t)}, & \alpha\in(1/4,\frac13], \\
            (1-\alpha)/2, & \alpha\in(\frac13,1),        
        \end{cases}
        \end{equation*}
        for $\alpha\in(1/4,1)$, Theorem~\ref{thm:main} immediately follows from Theorem~\ref{thm:main_gen}. Theorem~\ref{thm:main_crit} follows from Theorem~\ref{thm:main_gen} by choosing $\eta_0=\infty$. 
        \end{proof}
        
\appendix
    \section{Many-to-one and many-to-two lemmas}
    \begin{lemma}[many-to-one lemma]\label{lem:many-to-one}
        Let $(p_k)_{k\in\N}$ fulfill \eqref{eq:mom_offspr} and let $\xi\,:\,\R\to (0,\infty)$ be a locally integrable function, then for any $t\geq0$ and measurable non-negative function $f$, we have
        \begin{equation*}
            \Ee{x}{\xi}{\sum_{\nu\in N(t)} f(X^\nu_r\,;\,r\leq t) } = \E{x}{{\rm e}^{\int_0^t \xi(B_r)\,\mathrm{d}r}f(B_r\,;\,r\leq t)},
        \end{equation*}
        where $B$ is a Brownian motion.
    \end{lemma}
    \begin{lemma}[many-to-two lemma]\label{lem:many-to-two}
        Let $(p_k)_{k\in\N}$ fulfill \eqref{eq:mom_offspr} and recall the definition of $m_2$. Furthermore let $\xi\,:\,\R\to (0,\infty)$ be a locally integrable function, then for any $t\geq0$ and measurable non-negative function $f$, we have
        \begin{align*}
            &\Ee{x}{\xi}{\sum_{\substack{\nu_1,\nu_2\in N(t) \\ \nu_1\neq\nu_2}} f(X^{\nu_1}_r\,;\,r\leq t) f(X^{\nu_2}_r\,;\,r\leq t) }\\
            &= (m_2-2) \int_0^t \E{x}{{\rm e}^{\int_0^t\xi(B_r^{1,s})\,\mathrm{d}r + \int_s^t\xi(B_r^{2,s})\,\mathrm{d}r} \xi(B^{1,s}_s)f(B_r^{1,s}\,;\,r\leq t)f(B_r^{2,s}\,;\,r\leq t) }\,\mathrm{d}s,
        \end{align*}
        where $B^{1,s}$ and $B^{2,s}$ are two Brownian motions that coincide until time $s$ and move independently afterwards.
    \end{lemma}
    For proofs, see \cite[Lemma~1]{harris_robert_multiple_spines}.

	\section{Hitting time}
    The following result can be found in the literature; since its proof is very short, we include it here for completeness.
    \begin{lemma}\label{lem:HittingTimeDensity}
		Let $H_0$ denote the first time a Brownian motion $B_t$ reaches zero. For any $s>0$ and every $z\in\R\setminus\{0\}$ we have the explicit identity
		\begin{equation*}
			\Pb{z}{H_0 \in \mathrm{d}s} = \frac{1}{\sqrt{2\pi}} \frac{\abs{z}}{s^{3/2}} \exp\left\{- \frac{z^2}{2s}\right\}
		\end{equation*}	
	\end{lemma}
	\begin{proof}
		Due to symmetry it is sufficient to consider the case $z>0$. For this purpose consider
		\begin{align*}
			\Pb{z}{H_0 \leq s} &= \Pb{0}{\inf_{0\leq r\leq s} B_r \leq -z} = 2 \Pb{0}{B_{s}\leq -z} = 2 \Pb{0}{B_1 \leq -\frac{z}{\sqrt{s}}}
		\end{align*}
		by the reflection principle. Therefore,
		\begin{align*}
			\frac{\mathrm{d}}{\mathrm{d} s} \Pb{z}{H_0 \leq s} &= \frac{2}{\sqrt{2\pi}} \exp\left\{ -\frac{z^2}{2s} \right\} \frac{1}{2}  \frac{z}{s^{3/2}} = \frac{1}{\sqrt{2\pi}} \frac{z}{s^{3/2}} \exp\left\{ -\frac{z^2}{2s} \right\}.
		\end{align*}
	\end{proof}
	
	\section{A reflection principle for BBM}
	\begin{lemma}\label{lem:BBM_refl}
		Let $\xi \, : \, \R\to(0,\infty)$ be a locally integrable function and $y<x$. Then
		\begin{equation}\label{eq:BBM_refl}
			\Pbe{x}{\xi}{M_t \leq y} \leq \Pbe{x}{\xi}{\underline{M}_t \leq y} \leq 2\Pbe{x}{\xi}{M_t \leq y},
		\end{equation}
		where $\underline{M}_t = \inf_{s\in[0,t]} M_s$ denotes the running minimum of the leftmost particle position.
	\end{lemma}
	\begin{proof}
		The first inequality in \eqref{eq:BBM_refl} is apparent. For the second, observe that
		\begin{align*}
			\Pbe{x}{\xi}{\underline{M}_t \leq y} = \Pbe{x}{\xi}{M_t \leq y} + \Pbe{x}{\xi}{\underline{M}_t \leq y, M_t > y}.
		\end{align*}
		Let $h_y = \inf \{ s>0 \,:\, M_s = y \}$ be the first hitting time of $y$ and let $\mathcal{F}_{h_y}$ be the associated stopping-time $\sigma$-algebra. By the strong Markov property,
		\begin{align*}
			\Pbe{x}{\xi}{\underline{M}_t \leq y, M_t > y} = \Ee{x}{\xi}{\1{\underline{M}_t \leq y} \Pbe{x}{\xi}{ M_t > y \mid \mathcal{F}_{h_y} }}.
		\end{align*}
        Conditionally on $\mathcal{F}_{h_y}$, there is a particle at position $y$ at time $h_y$. Hence
        \[
            \Pbe{x}{\xi}{ M_t > y \mid \mathcal{F}_{h_y} } \leq \Pb{y}{ B_{t-s} > y}_{s={h_y}} = 1/2,
        \]
		and therefore
        \[
            \Pbe{x}{\xi}{\underline{M}_t \leq y, M_t > y} \leq \frac{1}{2} \Pbe{x}{\xi}{\underline{M}_t \leq y}.
        \]
        Substituting this bound into the previous decomposition yields
		\begin{align*}
			\Pbe{x}{\xi}{\underline{M}_t \leq y} \leq \Pbe{x}{\xi}{M_t \leq y} + \frac{1}{2} \Pbe{x}{\xi}{\underline{M}_t \leq y}
		\end{align*}
		and the proof is finished.
	\end{proof}

\bibliographystyle{alpha}
\bibliography{bibliography}

\begin{thebibliography}{NRRZ12}

\bibitem[AT07]{AdTa-07}
Robert~J. Adler and Jonathan~E. Taylor.
\newblock {\em Random fields and geometry}.
\newblock Springer Monographs in Mathematics. Springer, New York, 2007.

\bibitem[BH20]{BoHa-18}
Anton Bovier and Lisa Hartung.
\newblock From 1 to 6: a finer analysis of perturbed branching {B}rownian
  motion.
\newblock {\em Comm. Pure Appl. Math.}, 73(7):1490--1525, 2020.

\bibitem[Bra83]{bramson1983convergence}
Maury Bramson.
\newblock Convergence of solutions of the {K}olmogorov equation to travelling
  waves.
\newblock {\em Mem. Amer. Math. Soc.}, 44(285):iv+190, 1983.

\bibitem[BS02]{handbookBM}
A.~N. Borodin and Paavo Salminen.
\newblock {\em Handbook of {Brownian} motion: {Facts} and formulae.}
\newblock Probab. Appl. Basel: Birkh{\"a}user, 2nd ed. edition, 2002.

\bibitem[{\v{C}}DS23]{cerny2021unbounded}
Ji{\v{r}}{\'{\i}} {\v{C}}ern{\'y}, Alexander Drewitz, and Lars Schmitz.
\newblock (un-)bounded transition fronts for the parabolic {Anderson} model and
  the randomized {F}-{KPP} equation.
\newblock {\em Ann. Appl. Probab.}, 33(3):2342--2373, 2023.

\bibitem[Cha25]{LC25}
Louis Chataignier.
\newblock Upper moderate deviation probabilities for the maximum of branching
  {Brownian} motion.
\newblock {\em Electron. Commun. Probab.}, 30:14, 2025.
\newblock Id/No 87.

\bibitem[DS22]{DreSchmi19}
Alexander Drewitz and Lars Schmitz.
\newblock Invariance principles and log-distance of {F}-{KPP} fronts in a
  random medium.
\newblock {\em Arch. Ration. Mech. Anal.}, 246(2-3):877--955, 2022.

\bibitem[Fis37]{fisher1937}
R.~A. Fisher.
\newblock The wave of advance of advantageous genes.
\newblock {\em Annals of Eugenics}, 7:355--369, 1937.

\bibitem[HN22]{hamel2021diameters}
Fran\c{c}ois Hamel and Gr\'{e}goire Nadin.
\newblock Diameters of the level sets for reaction-diffusion equations in
  nonperiodic slowly varying media.
\newblock {\em Proc. Amer. Math. Soc.}, 150(8):3549--3564, 2022.

\bibitem[HR17]{harris_robert_multiple_spines}
Simon~C. Harris and Matthew~I. Roberts.
\newblock The many-to-few lemma and multiple spines.
\newblock {\em Ann. Inst. Henri Poincar\'e Probab. Stat.}, 53(1):226--242,
  2017.

\bibitem[INW68]{branch_markov_I}
Nobuyuki Ikeda, Masao Nagasawa, and Shinzo Watanabe.
\newblock Branching {M}arkov processes. {I}.
\newblock {\em J. Math. Kyoto Univ.}, 8:233--278, 1968.

\bibitem[KPP37]{kolmogorov1937}
A.~Kolmogorov, I.~Petrovskii, and N.~Piskunov.
\newblock Study of a diffusion equation that is related to the growth of a
  quality of matter and its application to a biological problem.
\newblock {\em Moscow University Mathematics Bulletin}, 1:1--26, 1937.

\bibitem[McK75]{mckean}
H.~P. McKean.
\newblock Application of {B}rownian motion to the equation of
  {K}olmogorov-{P}etrovskii-{P}iskunov.
\newblock {\em Comm. Pure Appl. Math.}, 28(3):323--331, 1975.

\bibitem[NRRZ12]{nolen_zlatos_12}
James Nolen, Jean-Michel Roquejoffre, Lenya Ryzhik, and Andrej Zlato\v{s}.
\newblock Existence and non-existence of {F}isher-{KPP} transition fronts.
\newblock {\em Arch. Ration. Mech. Anal.}, 203(1):217--246, 2012.

\bibitem[Pit99]{Pi-99}
Jim Pitman.
\newblock The distribution of local times of a {B}rownian bridge.
\newblock In {\em S\'eminaire de {P}robabilit\'es, {XXXIII}}, volume 1709 of
  {\em Lecture Notes in Math.}, pages 388--394. Springer, Berlin, 1999.

\bibitem[Roc24]{MR5088005}
S\'ebastien Roch.
\newblock {\em Modern discrete probability---an essential toolkit}, volume~55
  of {\em Cambridge Series in Statistical and Probabilistic Mathematics}.
\newblock Cambridge University Press, Cambridge, 2024.

\bibitem[Uch77]{uchiyama77}
K\^{o}hei Uchiyama.
\newblock The behavior of solutions of the equation of
  {K}olmogorov-{P}etrovsky-{P}iskunov.
\newblock {\em Proc. Japan Acad. Ser. A Math. Sci.}, 53(7):225--228, 1977.

\bibitem[Zla12]{zlatos12}
Andrej Zlato\v{s}.
\newblock Transition fronts in inhomogeneous {F}isher-{KPP} reaction-diffusion
  equations.
\newblock {\em J. Math. Pures Appl. (9)}, 98(1):89--102, 2012.

\end{thebibliography}

\end{document}